\documentclass[twoside,english]{amsart}
\usepackage{color}
\usepackage[utf8]{inputenc}
\usepackage[T1]{fontenc}
\usepackage{lmodern}
\usepackage{babel}
\usepackage{geometry} 
\usepackage{graphicx}

\usepackage{amsbsy}
\usepackage{amssymb,latexsym}
\usepackage{amsmath,amsthm}
\usepackage{mathtools}

\usepackage{pgffor}
\usepackage{tikz}
\usepackage{tikz-cd}
\usetikzlibrary{matrix,arrows}
\usetikzlibrary{intersections}
\newcommand{\C}{\mathbb{C}}
\newcommand{\Pp}{\mathbb{P}}

\newcommand{\M}{\overline{\mathcal{M}}}
\newtheorem{thm}{Theorem}

\newtheorem{defn}{Definition}

\newtheorem{lem}{Lemma}
\newtheorem{prop}{Proposition}
\newtheorem{cor}{Corollary}

\newtheorem{rem}{Remark}
\numberwithin{thm}{subsection}

\numberwithin{prop}{subsection}
\numberwithin{lem}{subsection}
\numberwithin{cor}{subsection}

\numberwithin{defn}{subsection}

\numberwithin{prop}{subsection}
\numberwithin{lem}{subsection}
\numberwithin{cor}{subsection}

\newtheorem{example}{Example}

 \def\dpol{\mathop{{}^{\textsc{D}}\text{Pol}}}

 \def\eps{{\varepsilon}}

\def\ur{{\underline r}}
\def\uv{{\underline v}}
\def\ueps{{\underline \varepsilon}}
 \def\dpol{\mathop{{}^{\textsc{D}}\text{Pol}}}
 \def\C{\mathbb{C}}
   \def\M{\overline{\mathcal{M}}}
  \def\Rl{\mathbb{R}}

\title{The Gauss--skizze decomposition is a Goresky--MacPherson stratification.}

\author{N.C. Combe}
\address{University of Warsaw}
\email{n.combe@uw.edu.pl}
\thanks{I wish to thank Leila Schneps, Pierre Lochak and Hans Henrik Rugh for helpful suggestions and discussions during my stay at Paris 6. All my gratitude goes towards Yuri Ivanovitch Manin who introduced me to moduli spaces $\M_{g,n}$ and Frobenius manifolds. I would like as well to thank very much the unknown referee for very relevant comments.}

\begin{document}
\maketitle

\begin{abstract}
We consider a new stratification of the space of configurations of $n$ marked points on the complex plane. Recall that this space can be differently interpreted as the space $\dpol_{n}$ of degree $n>1$ complex, monic polynomials with distinct roots, the sum of which is 0. A stratum $A_{\sigma}$ is the set of polynomials having $P^{-1}(\mathbb{R}\cup\imath\mathbb{R})$ in the same isotopy class, relative to their asymptotic directions. We show that this stratification is a Goresky--MacPherson stratification and that from thickening strata a {\it good cover} in the sense of \v Cech can be constructed, allowing an explicit computation of the cohomology groups of this space. 
\end{abstract}

\setcounter{tocdepth}{1}
\tableofcontents

\

\section{Introduction}
The decomposition of the space of configurations of $n$ marked points on the complex plane has been considered for over more than fifty years~\cite{FoN62} and has lead to many important works~\cite{{Ar69},{CMS08},{CPS01},{Coh73A},{Fu70},{Gor78},{FaN62},{Va78}}.
The braid group $\mathcal{B}_n$ with $n$ strands and the space $\dpol_{n}$  of degree $n$ complex monic polynomials with distinct roots are objects which are deeply connected: 
the space $\dpol_{n}$ is a  $K (\pi, 1)$ of fundamental group $\mathcal{B}_n$.
The richness of their interactions allows each object to provide information about the other.  
In particular, one may use the space $\dpol_{n}$ to give a description of braids with $n$ strands. By means of a spectral sequence V. Arnold~\cite{Ar69} gives a  method to calculate the integral cohomology of the braid group. Another version was given by D.B. Fuks~\cite{Fu70}, which allows one to calculate the cohomology of the braid group with values in $\mathbb{Z}_{2}$. This latter approach was further on developed in a more general way by F.V.  Vainshtein~\cite{Va78} using Bockstein homomorphisms. More recent results on the topic have appeared in the works of V. Lin, F. Cohen~\cite{Coh73A} and  A. Goryunov~\cite{Gor78}. For results more closely related to ours, we also mention the work of C. De Concini, D. Moroni, C. Procesi, M. Salvetti~\cite{CPS01,CMS08}. 

\smallskip
In spite of an abundant literature concerning this subject, we give a new approach to the configuration space of $n$ marked points on $\C$ and thus on the moduli space $\M_{0,n}$ of $n$ marked points on the Riemann sphere $\mathbb{P}^1$. 
In this paper, we develop the {\it real}-geometry of the canonical stratification of the space $\M_{0,n}$, which gives a new insight on \cite{CoMa,De,Ke,KoMa,Ma99}. 

\smallskip 

To be more precise, the aim is to present a real-algebraic stratification of the parametrizing  $\C$–scheme of $\M_{0,n}$. We show that it is a {\it Goresky--MacPherson stratification}.
This highlights real-geometry properties of $\M_{0,n}$ and thus gives a different approach to this object. 

\smallskip 

The notion of Goresky--MacPherson stratification  can be found in \cite{GM0,GM1,GM2}. 
Let us recall the notion of stratification. 

\smallskip

A stratification of a topological space $X$ is a locally finite partition of this space $\mathcal{S}=(X^{(\sigma)})_{\sigma\in S}$ of $X$ into elements called strata, locally closed and verifying the condition: for all $\tau,\sigma\in S,$ 

\[X^{(\tau)}\cap \overline{X}^{(\sigma)}\neq \emptyset \iff X^{(r)}\subset \overline{X}^{(\sigma)}.\]

The boundary condition can be reformulated as follows: the closure of a stratum is a union of strata.
A different interpretation of this would be to define a stratification as a filtration
  \[X=X_{(n)}\supset X_{(n-1)}\supset \dots \supset X_{(1)}\supset X_{(0)}\supset X_{-1}=\emptyset,\]
where the $X_{(i)}$ are closed, and where by defining $X^{(i)}:=X_{(i)}\setminus X_{(i-1)},$ we have $\overline{X^{(i)}}=X_{(i)}$.

\smallskip 

The Goresky--MacPherson stratification is defined as follows. An $n$-dimensional topological stratification of a paracompact Hausdorff topological space $X$ is a filtration by closed subsets
\[X=X_{(n)}\supset X_{(n-1)}\supset \dots \supset X_{(1)}\supset X_{(0)}\supset X_{-1}=\emptyset,\]
such that for each point $p\in X_{(i)}-X_{(i-1)}$ there exists a distinguished neighborhood $N$ of $p$ in $X$, a compact Hausdorff space $L$ with an $n-i-1$ dimensional topological stratification 
\[L=L_{n-i-1}\supset \dots \supset L_1\supset L_0\supset L_{-1}=\emptyset\]
and a homeomorphism 
\[\phi:\Rl^i\times cone^\circ(L)\to N,\] which takes $\Rl^i \times cone^{\circ} (L_j) \to N\cap X_{i+j+1}$. The symbol  $cone^{\circ} (L_j)$ denotes the open cone, $L\times [0,1)/(l,0)\sim(l',0)$ for all $l, l'\in L$.

\smallskip 

Focusing on the case where the marked points on $\C$ are pairwise distinct, we show that from this stratification can be constructed a {\it good \v Cech cover.} Moreover, this paper provides an explicit construction of the open sets of the \v Cech cover. 

\smallskip

\begin{defn}[\v Cech cover~\cite{Cech32}]
Consider a cover of a topological space $X$. A good \v Cech cover is a cover such that its components are open and have contractible multiple intersections.\end{defn}  

\smallskip

We prove the following theorems. 
\begin{thm}
Consider the configuration space of $n$ marked points on the complex plane, where points are pairwise disjoint. 
Then, there exists a real algebraic stratification $\mathcal{S}$ of this space, where strata are indexed by oriented and bicolored forests verifying the following properties:
\begin{itemize}
\item there exist $n$ vertices of valency 4, incident to edges of alternating color and orientation;
\item there exist at most $n-1$ vertices of even valency, incident to edges of only one color and of alternating orientations;
\item there are $4n$ leaves (i.e. vertices incident to one edge), where the colors and orientations of the incident edges alternate. 
\end{itemize}
\end{thm}

\begin{thm}
Let $A_\sigma$ be a generic stratum (i.e. of codimension 0). Then, the topological closure $\overline{A}_\sigma$ defines a Goresky--MacPherson stratification.
\end{thm}

Thickening these strata prepares the ground for defining a good \v Cech cover, which leads to the following statement:

\begin{thm}
Consider the configuration space of $n$ marked points on the complex plane, where points are pairwise disjoint, and the stratification $\mathcal{S}$ defined above.
Then, the thickened strata $\underline{A_{\sigma}}^+$ form a good \v Cech cover. 
\end{thm}

\smallskip 

In a more global way, what we have in mind is a new manner of defining the generators for $\Gamma_{0,[n]}$, where $\Gamma_{0,[n]}$ is nothing but the orbifold fundamental group of the moduli space of smooth curves of genus  0 with $n$ unordered marked points~\cite{Gro84,LS}. 
This new cell decomposition turns out to have many interesting symmetries (this is the subject of the paper~\cite{CO2}), where we aim at underlining the existence of polyhedral symmetries explicitly in the presentation. In this paper we show the existence of a good cover in the sense of \v Cech which we will use, in further investigations, to pave the way {  towards persistent homology}.
 
 \smallskip
 
Namely, to stratify the configuration space $\dpol_n$, we consider {\it drawings of polynomials}---objects reminiscent of Grothendieck's {\it dessin's d'enfant} in the sense that we consider the inverse image of the real and imaginary axis under a complex polynomial. In particular, since Gauss was first to discover these special graphs, we call them {\it Gauss skizze} which means Gauss' drawings.

The {\rm drawing} associated to a complex polynomial is a system of blue and red curves properly embedded in the complex plane, being the inverse image under a polynomial $P$ of the union of the real axis (colored red) and the imaginary axis (colored blue)~\cite{CO1}. For a polynomial of degree $d$, the
drawing contains $d$ blue and $d$ red curves, each blue curve intersecting exactly one red curve
exactly once.  The entire drawing forms a forest whose leaves (terminal edges) go to infinity
in the asymptotic directions of angle $k \pi/4d$.  Polynomials belong to a given stratum if their drawings are isotopic, relatively to the $4d$ asymptotic directions.  An element of the stratification is indexed by an equivalence class of drawings which we call a signature $\sigma$. A stratum shall be denoted by $A_{\sigma}$. Each stratum of this decomposition is attributed to a decorated graph  (i.e. a graph where edges are oriented and colored in red or blue; complementary faces to the embedded graph are colored in colors $A,B,C,D$). 

\smallskip 
{ 
For clarity such a decorated forest is visualized as a chord diagram in the unit disc. It is obtained by an embedding of the forest to a closed unit disc such that leaves are mapped to the boundary circle and the rest of the forest is mapped into the open disc. One places the leaves at the $4n$-roots of unity on the unit circle.}

\smallskip

 This paper is organized as follows.  Section 2, we present the definition and construction of the real-algebraic stratification. We introduce the notion of drawings, signatures and present some of its properties. 

\smallskip

Section 3 introduces Whitehead moves on the signatures which allow a definition of the combinatorial closure of a signature. An important result which follows from this section is that the topological closure of a stratum is given by its combinatorial closure. In particular, it is the union of all the incident signatures to $\sigma$ ({\it i.e} the combinatorial closure). Finally, we prove that this stratification is a Goresky--MacPherson stratification.

\smallskip

Section 4 shows the construction of the \v Cech cover. In particular, we define the components of the cover and prove that multiple intersections between those components are contractible. 

Finally, an appendix is introduced, where we discuss the multiple intersections of closures of generic strata and classify the possible diagrams in order to have non-empty multiple intersections.

\medskip

\section{A new point of view on $\M_{0,T}$}
\subsection{Moduli spaces of genus 0 curves with marked points}
Let $\C$ be  fixed field and $\bf{T}$ a finite set. Consider the moduli space (generally a stack) $\M_{0,\bf{T}}$ of stable genus zero curves with a finite set of smooth pairwise distinct closed points defined over $\C$ and bijectively marked (labeled) by $\bf{T}$. For a $\C$–scheme $B$, a $B$–family of such curves is represented by a fibration $C_B \to B$ with genus zero fibres of dimension 1, endowed by pairwise distinct sections $s_i : B \to C$ labeled by $\bf{T}$ as above (cf. \cite{Ma99}).

\smallskip 

From the known complex results, the combinatorial type of each such curve over a closed point is encoded by a stable tree $\sigma$. Leaves (or else, tails) of such a tree are labeled by the elements of a subset $\bf{T}_\sigma \subseteq \bf{T}$. It has been considered in \cite{KoMa} and later in \cite{Ma99} and \cite{CoMa} the spaces $\mathcal{M}_{0,\boldsymbol{\sigma}}$ and their natural closures $M_{0,\boldsymbol{\sigma}} \subset \M_{0,\boldsymbol{\sigma}}.$

\smallskip

One point of $\M_{0,n+1}(\C)$ “is” a stable complex curve of genus 0 with $n + 1$ points labeled by $\{1,2,...,n-2, 0, 1, \infty\}$. These points are distributed along different locally closed strata that are naturally labeled by stable trees. The complex codimension of a stratum is the number of edges of the respective tree. We will refer to this number as level. Each stratum is reduced and irreducible, two different strata do not intersect.

\smallskip 

At the level zero, all these pairwise different points live on a fixed $\Pp^1$ endowed with a fixed real structure, with respect to which $(0,1,\infty)$ are real. Equivalently, this projective line is endowed with a fixed homogeneous coordinate system $(y : z)$ such that $0=(0:1)$, $1=(1:1)$, $\infty=(1:0)$.

\smallskip 

\subsection{Definition-construction}

Given a stratification of $\mathcal{S}=(X^{(\sigma)})_{\sigma \in S}$ one can define an incidence relation on the set $S$ of indices by: 
\[\tau \dot{\preceq} \sigma\quad  \text{if and only if} \quad X_{(\tau)}\cap\overline{X}^{(\sigma)}\neq \emptyset.\]
In other words, in a stratification the incidence relation gives a {\it poset}.

\smallskip

\paragraph{We consider the set of $n+1$ marked points on $\Pp^1$ as a configuration space of $n$ points on the complex plane modulo the group $PSL_2(\C)$. This configuration space $Conf_{n}(\C)$ can itself be considered as the set of monic complex polynomials, with $n$ roots, which we denote $\dpol_{n}$. The quotient of this space by $PSL_2(\C)$ allows to map three of the marked points to $\{0,1,\infty\}$.} 

\smallskip

\paragraph{For the construction of this new stratification, it is very important to have a fixed real structure on $\Pp^1$. The main idea of our construction is to take the inverse image, under a polynomial $P$ in $\dpol_n$, of the real and imaginary axis, i.e. $P^{-1}(\Rl\cup\imath \Rl)$. This inverse image forms a system of oriented curves in the complex plane. To distinguish $P^{-1}(\Rl)$ from $P^{-1}(\Rl\imath\Rl)$, we color in red the pre-image of $\Rl$ and in blue the pre-image of $\Rl\imath\Rl$. The orientation of the curves is inherited from the natural orientation on the real and imaginary axis (see Figure \ref{F:patdisc}).}

\smallskip
\begin{defn}
A (Gauss) $n$-drawing $\mathcal{C}_{P}$  of a degree $n$-polynomial $P\in \dpol_n$ is a system of curves properly embedded in the complex plane given by $P^{-1}(\Rl\cup\imath \Rl)$.
\end{defn}
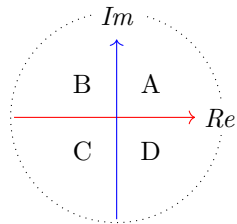
\begin{figure}[h]
\begin{center}
\begin{tikzpicture}[scale=0.45],
every node/.style = {inner xsep=0mm, inner ysep=1mm, fill=white}
\draw[dotted] (0,0) circle (3.1cm);
\draw [->,red] (-3,0) -- (2.3,0);
\draw (3,0) node[fill=white]{\it Re};
\draw[->,blue] (0,-3) -- (0,2.3);
\draw (0,3) node[fill=white]{\it Im};
\draw (1,1) node[fill=white]{A};
\draw (-1,1) node[fill=white]{B};
\draw (-1,-1) node[fill=white]{C};
\draw (1,-1) node[fill=white]{D};
\end{tikzpicture}
\vspace{-3pt}
\caption{Partition of the complex plane}\label{F:patdisc}
\end{center}
\end{figure}
\smallskip
For simplicity we will just refer to Gauss' drawings as drawings. 
\begin{rem}\label{R:1}{ An important remark concerning the topology of the drawings follows from the tight relations between a complex polynomial and harmonic polynomials.  Indeed, any degree $n$ complex polynomial $P$ splits into a sum of harmonic polynomials in two variables $\Re(P)$and $\Im(P)$ given by the formula $P=\Re(P)+\imath \Im(P)$, where $\Re(P)$ is the real part of $P$ and $\Im(P)$ is the imaginary part of $P$. The level sets of a given harmonic polynomial $u$ are given by
$\{z:\, u(z)=0\}$ and forms an embedded forest with $2n$ leaves and inner nodes of even valencies. The non-existence of ovals follows from the fact that the maximum principle would be otherwise contradicted. }
\end{rem}

{  We regard the drawings as embedded graphs (forests) in the complex plane. Those are forests with $4n$
 leaves lying on the the roots of unity. In other words we have $4n$ asymptotic directions lying on the roots of unity.}

Using the coloring convention of the curves above, we now list the different types of intersections of curves in a drawing:

\smallskip 

\begin{enumerate}

\item The {\it roots} of a polynomial. These are given by the intersection of a red curve with a blue curve.  \\

\item The {\it critical points} of a polynomial.
\begin{itemize}
\item Let $z_0$ be a critical point of $P$ such that its corresponding critical value $P(z_0)\in \imath \Rl$ is an imaginary number. This is satisfied for the algebraic equation $\Re(P)(z_0)=0$. The point $z_0$ lies at the intersection of the blue curves. \\
\item Let $z_0$ be a critical point of $P$ such that the critical value associated to $P(z_0)\in \Rl$ is a real number. This is satisfied for the algebraic equation $\Im(P)(z_0)=0$. The point $z_0$ lies at the intersection of the red curves.\\
\item  Let $z_0$ be a root of multiplicity $k$ i.e one has $P(z)=(z-z_0)^kg(z)$, where $g(z)$ is a polynomial. Then, the point $z_0$ is the intersection point of $k$ red and $k$ blue curves. The colors of the incident curves to $z_0$ alternate: red, blue, red, blue (etc); whereas the number of incident incoming and outgoing (half) curves is equal to $4k$. 
\end{itemize}
\end{enumerate}
The last case arrises only when one considers the discriminant variety $\Delta=\{z_i= z_j\, |\, z_i\in \C\}$; the two cases listed above it arise in the complement of the discriminant variety.  
\smallskip 
{ 
Those drawings are properly embedded forests with $4n$ leaves lying on the the roots of unity, that is we have $4n$ asymptotic directions lying on the roots of unity. We are sort of compactifying $\C$ with the circle of directions at infinity. Moreover, the rays themselves need not be asymptotic in the sense of getting closer together at infinity, merely that they stay bounded distance to a line of a fixed slope. Those asymptotic directions are preserved under homeomorphism, implying that it is a proper isotopy. This allows to define the following equivalence class of drawings.}
\begin{defn}\label{D:equi}
Two $n$-drawings $\mathcal{C}_{P_{1}}$ and $\mathcal{C}_{P_{2}}$ are equivalent if and only if there exists a proper isotopy $h$ of $\mathbb{C}$ (a continuous family of homeomorphisms of $\mathbb{C}$) such that $h:\mathcal{C}_{P_{1}}\to\mathcal{C}_{P_{2}}$ and $h$ preserves the $4n$ asymptotic directions, the colouring and orientation of curves. 
\end{defn} 

We denote the equivalence class of isotopic drawings by $[\mathcal{C}_{P}]$. This definition serves further on in the construction of the decomposition of the configuration space of marked points on the complex plane. 

\begin{defn} 
Consider the space of configurations of $n$ marked points on $\C$.
Let $A_{[\mathcal{C}_{P}]}$  be the set of polynomials with drawings in the isotopy class $[\mathcal{C}_{P}]$.  The family $(A_{[\mathcal{C}_{P}]})_{[\mathcal{C}_{P}] \in S}$, where $S$ is the set of isotopy classes of $n$-drawings, partitions the configuration space. The component $A_{[\mathcal{C}_{P}]}$ is called stratum.
\end{defn} 

\begin{defn}[Codimension]\label{D:cod}
Let $P$ be a degree $n$ complex polynomial in $\dpol_n$. 
\begin{enumerate}
\item A polynomial $P$ with no critical values on $\mathbb{R}\cup \imath \mathbb{R}$ is called generic, such a polynomial is of codimension 0.
\item The {\rm special} critical points of $P$ are the critical points $z$ such that $P(z) \in \mathbb{R}\cup \imath \mathbb{R}$. 
The local index at a special critical point $z\in \mathbb{R}$ (resp. $ \imath \mathbb{R}$) is equal to $2m-3$, where $m$ is the number of red (or blue) diagonals crossing at the point $z$. The {\rm  real codimension} of $P$ is the sum of the local indices of all the special critical points. 
\end{enumerate}
\end{defn}
\smallskip
\begin{example}
Figure~\ref{SSSS6} illustrates a generic (codimension 0) polynomial of degree 6.
\begin{figure}[h]
\begin{center}
\includegraphics[scale=0.3]{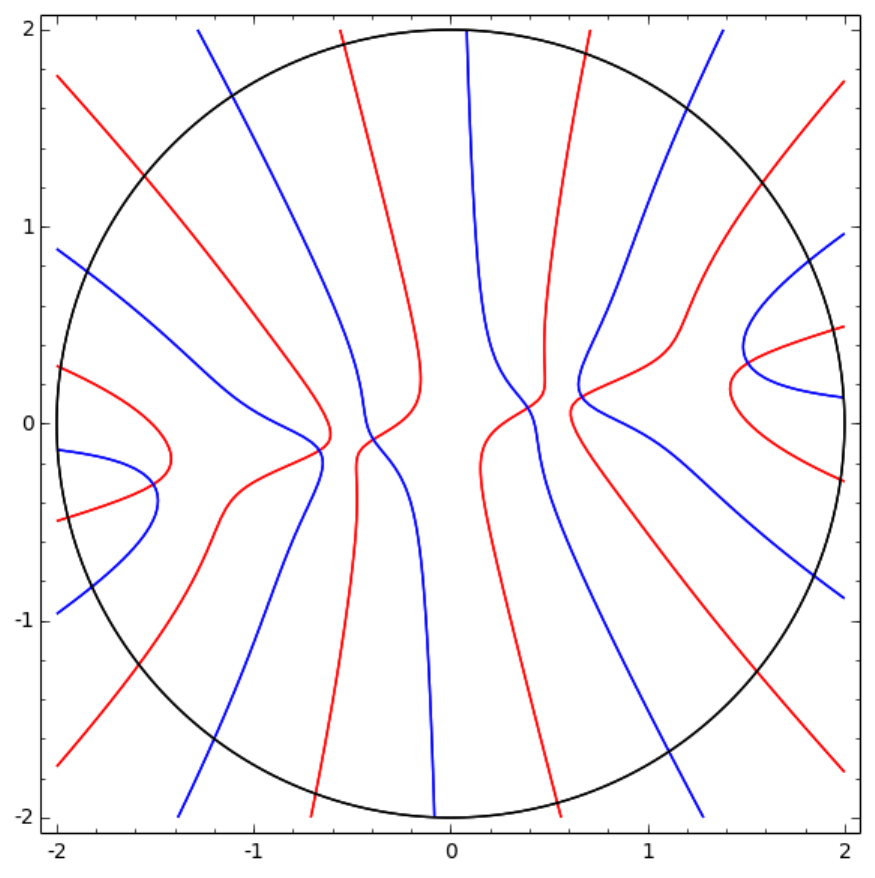} 
\end{center}
\vspace{-10pt}
\caption{$P(z)= (z+1.54*(0.98+i*0.2))*(z+0.68*(0.98+i*0.2))*(z+0.4*(0.98+i*0.2))*(z-1.54*(0.98+i*0.2))*(z-0.68*(0.98+i*0.2))*(z-0.4*(0.98+i*0.2)).$}\label{SSSS6}
\end{figure}
\end{example}
\medskip

\subsection{Embedded forests}

Consider the category $\Gamma$ of graphs. Objects of this category are graphs, i.e. objects defined by a set $V$ of elements called vertices, equipped with a  symmetric, reflexive, relation 
$E$. For given vertices $x,y\in V$ there is an edge $(x,y)$ if and only if $(x,y)\in E$ and $x\ne y$.  We define a morphism $(V,E)\to (W,F)$ as a function $f: V\to W$ such that $(x,y)\in E$ implies $(f(x),f(y))\in F$. 

\smallskip 

We equip this category $\Gamma$ with a coproduct $\sqcup$, being here the disjoint union of graphs. In the category $\Gamma$, we are interested only in those graphs being trees, i.e. acyclic connected graphs. By convention, a tree has at least one edge. Vertices of valency one are called {\it leaves}. 
We equip the standard definition of trees with a colouring and an orientation on the set of edges. Each edge is mapped to the set $\{R^{+},B^{-},B^{+},R^{-}\}$, where the capital letter $R$ (resp. $B$) corresponds to a red (resp. blue) colouring and the sign "+" or "-" corresponds to the orientation of the edge.  We will write $R^{-}$ (or $B^{-}$) if at a given vertex, the direction points inwards, i.e. towards the chosen vertex. Similarly, if we have the opposite sign, then the direction points in the opposite direction i.e. outwards.
Concerning the notation, we use the symbol $(i,j)_{R}$ (resp. $(i,j)_{B}$) for a red (blue) diagonal in the forest, connecting the leaf $i$ to the leaf $j$. 

We define the following object. 
\begin{defn}
An $n$-signature is an object of the category $(\Gamma,\sqcup)$ being a disjoint union of trees (i.e. a forest)
such that: 
\begin{enumerate}
\item There exist $k$ vertices of valency four. The incident edges are in bijection with the set $\{R^{+},B^{-},B^{+},R^{-}\}$. Orientation and colors alternate. 
\item There exist $n-k$ vertices of valency a multiple of four. The incident edges are in bijection with the set $\{R^{+},B^{-},B^{+},R^{-}\}$. Orientation and colors alternate. 
\item There exist at most $k-1$ vertices of even valency (at least 4), incident to edges of the same color. Orientations alternate.
\item There exist $4n$ leaves. The set of incident edges to those leaves are of alternating color and orientation.  
\end{enumerate}
\end{defn}
Note that for the case of pairwise distinct points on the complex plane we have exactly $n$ vertices incident to four edges of alternating colors/orientations and at most $n-1$ vertices incident to monochromatic edges.

 \medskip
 
\begin{itemize} 
\item A\ {\it geometric realisation} of a  signature is a (union of) 1-dimensional CW-complex, being contractible, and properly embedded in the complex plane. 
\item An {\it embedded forest} is a subset of the plane which is the image of a proper embedding of a forest minus leaves to the plane. 
\end{itemize}

It is important for the construction presented in this note, to insist on the difference between the geometric realisation of a given graph and the graph itself, which is a purely {\it combinatorial} notion. In particular, for clarity we shall visualise a signature as a chord diagram i.e. a graph embedded to a closed unit disc such that the leaves are mapped to the boundary circle and the rest of the forest into the open disc. We place the leaves at the $4n$-roots of unity on the unit circle.

\medskip

\subsection{From real curves to forests}
 
An equivalent way of describing the orientation on the edges on the graphs is to consider the complementary part to the real and imaginary axis in $\C$ and to color them. {  We call the connected components of the complement of the embedded forest a face. The degree of a face is defined as the number of connected components of its boundary.}Let us label the four quadrants in $\C\setminus\Rl\cup\imath\Rl$ in the colors $A,B,C,D$ as in the figure \ref{F:patdisc}.

Relating to the previous remark~\ref{R:1} on the relation between harmonic polynomials and their level curves being forests we have the following lemma. 
\begin{lem}
Let  $\gamma\in [\mathcal{C}_{P}]$, then $\gamma$ is a forest.
\end{lem}
\begin{proof}
Consider an embedded forest $F_R$ (resp. $F_B$) with $2n$ leaves, whose edges are colored in $R$ (resp. B) and such that all vertices in the plane are of even valencies. Then, according to theorem 1.3  in \cite{Er}, there exists a harmonic polynomial of degree $n$ whose zero set is equivalent to $F_R$ (resp. $F_B$). The class $[\mathcal{C}_{P}]$ is a forest, since it is easy to see that the level set of the degree $n$ harmonic polynomials $\{(x,y)\in \Rl^2: \Re P(x,y)=0\}$ (resp. $\{(x,y)\in \Rl^2: \Im P(x,y)=0\}$) is an embedded forest (the existence of a cycle would contradict the maximum principle).  
\end{proof}

\begin{thm}\label{T:1}
The set of $n$-signatures is in bijection with the set of isotopy classes of drawings relatively to the $4n$ asymptotic directions.
\end{thm}
\begin{proof}
One direction is easy.  The other direction is proved as follows. Let $\gamma$ be the embedding in $\C$ of a signature $\sigma$ 
Then, we can construct a function $f$ such that:
\begin{itemize}
\item $f: \C \to \C$ is a smooth function,
\item $f^{-1}(\Rl\cup\imath\Rl)=\gamma$
\item The function $f$ is a bijection between each region of $\gamma$ and a quadrant colored $A, B, C, D$
in the complex plane,
\item $f$ is injective and regular on the edges of $\gamma$.
\end{itemize}
Let $J$ be the standard conformal structure on $\C$. We have $J= f_{*}(J_{0})$ so $f:(\C,J)\to (\C,J_{0})$ is holomorphic. By Riemann’s mapping theorem there exists a biholomorphic function $\rho: (\C, J ) \to (\C, J_0)$. The classical theorem of complex analysis by Rouch\'e implies that  $f \circ \rho^{-1}$ is a polynomial.
\end{proof}

\smallskip

We now proceed to the proof of the fact that a stratum $A_{\sigma}$ is contractible. The tools to prove it require the introduction of some statements due to Cerf~\cite{Ce62}. 

\begin{thm}\label{T:Contract}
Let $\sigma$ be an $n$-signature. Then, the set $A_{\sigma}$ is contractible.
\end{thm}

The theorem we need to use (and to adapt) is a theorem of J. Cerf. This statement applies in general to manifolds with boundary. 
Consider $V$ a manifold with boundary and compact.  
Let $W$ be a "target" manifold. Denote bt $Pl(V,W)$ the space of smooth embeddings of $V$ into $W$.
Then, we have the following. 
\begin{thm}[Cerf, theorem 3, \cite{Ce62}]\label{T:C} 
Let $G$ be a subgroup of the group of all diffeomorphisms of $V$, where $G$ acts on $Pl(V,W)$ and determines a structure of principal bundle. If $G$ is open, then the canonical map:
\[Pl(V,W)\to Pl(V,W)/G\] is a locally trivial fibration. 
\end{thm}
For an exposition on principal bundles, we refer to \cite{Ca}.

\smallskip 
{  We now comment on this construction and how we intend to make use of it. It is the graph $\sigma$ that plays the role of the manifold with boundary (one has a 1-dimensional manifold where the boundaries are given by the fixed intersection patterns between the red and blue curves) and which embeds in the target manifold $W=\C$. The smooth embeddings of the signatures $\sigma$ in $\C$ are denoted by $Pl(\sigma,\C)$. It is important to remark that this statement holds only for transversal intersections. In particular, if one considers two curves which meet at a tangency; they don’t intersect for nearby embeddings. This is exactly what happens when a critical point with real or pure imaginary critical value arises.
However, one can easily handle this situation by using the fact that for such critical points with real or pure imaginary critical values, one can interpret them as transverse intersections: when one has such a point there exists an even number of prongs. So, one can pair them off as if pairwise transverse intersections would be given.

}

Moreover, the following lemma (Lem. 1 in \cite{Ce62}) is used in our proof of the contractibility theorem and we recall it for the convenience of the reader. 

\begin{lem}\label{L:1}
Let $(E,B, p)$ be a fiber bundle where $E$ and $B$ are topological spaces, $p$ is a  continuous map $p: E \to B$. 
A sufficient condition to have a locally trivial fibration is that for $x_{0}\in B$ there exists a topological group $G$ operating (left wise) on the total space $E$
and on the base space $B$ so that the following diagram is commutative:

\[
\begin{tikzcd}
G \times E \arrow{r}{\phi} \arrow[swap]{d}{Id \times p} & E \arrow{d}{p} \\
G \times B  \arrow{r}{\Phi} & B
\end{tikzcd}
\]
and the map $g \to g x_{0}$ from $G$ into $B$ has a continuous section in the neighborhood of $x_{0}$.

\end{lem}

\begin{rem}\label{R:1}Let us comment on this result. 
Let $\gamma$ be the embedding in $\C$ of a signature $\sigma$, i.e. $\gamma \in Pl(\sigma,\C)$ which satisfies the asymptotic directions. 
Let $T$ be a tubular neighbourhood closed of $\gamma$ in $\C$ and denote by $H$ the group of diffeomorphisms, inducing the identity on $\C-T$.   

\smallskip
 From \ref{L:1}, it is shown in \cite{Ce62} that there exists a neighborhood $\mathcal{V}$ of $\gamma$ in $Pl(\sigma,\C)$ and a continuous map  $s: \mathcal{V} \to H$ such that :  \begin{itemize}
\item $s(\gamma)=e;$
\item For all $\gamma'\in \mathcal{V}$, $s(\gamma')\gamma'$ is identified with a diffeomorphism of $\sigma$ in the neighborhood of the identity, which is of the form $\gamma g$ where $g$ is a diffeomorphism of $\sigma$.
\item $\forall \gamma'\in \mathcal{V}$ and any diffeomorphism $g$ of $\sigma$ such that $\gamma' g\in \mathcal{V}$, $s(\gamma' g)=s(\gamma')$
\end{itemize} 

\end{rem}

\smallskip 

\begin{proof} 
We will mainly use the Cerf theorem and the remark \ref{R:1}, above. 

\smallskip

The embedded graphs of $\sigma$ in $\C$ correspond to a class of drawings. We know that this class is contractible, by  theorems 5.3 and 5.4 of Epstein in \cite{Ep}. 
For each $\gamma$ in the neighborhood $V$, we can construct the function $f$, as mentioned above in theorem \ref{T:1}. Define $E_{\gamma}$ to be the space of those functions corresponding to the embedded graph $\gamma$, and verifying the properties described in the proof of theorem \ref{T:1}. This space is contractible. 

Define the Cerf fibration \[E_\sigma=\{ (\gamma,f) | \gamma \in \sigma ,f \in E_{\gamma}\}\to \sigma.\]
\[(\gamma,f )\mapsto \gamma\]

\smallskip

By the result of Cerf in Theorem \ref{T:C}, we have a structure of principal bundle on $E_\sigma \to E_\sigma/G$, where $G$ is the group of oreintation-diffeomorphisms $G$ acting on $\C$, which are homotopic to the identity. Furthermore, since the isotopies need to preserve the circle at infinity  
by a theorem of Smale, this group is contractible \cite{Sma}.
Now, by a theorem of Cerf \cite{Ce61} the total space $E_{\sigma}$ is contractible. Note that $G$ acts without fixed points on $E$ and that orbits are closed. 
Using the Riemann mapping theorem and the classical Rouch\'e theorem, implies that each class modulo $G$ has a unique representative which is a complex polynomial.

\end{proof}

\section{From real curves to stable trees}

\subsection{Stratification on $\M_{0,n}$}
Let $\boldsymbol{\tau}$ be the (labeled) graph of a stable connected curve with $n$ marked points $C_B$. Denote by $M_{0, {\boldsymbol{\tau} }} \subset \M_{0,n}$ the moduli submanifold (or generally, substack)
parametrizing  all curves having the same graph ${\boldsymbol{\tau}}$.  Its closure will be denoted
$\M_{0,\boldsymbol{\tau}}$.

\smallskip

A morphism $f:\, {\boldsymbol{\tau}} \to {\boldsymbol{\sigma}}$ is determined by its covariant surjective action
 upon vertices $f_v:\, V_{{\boldsymbol{\tau}}}\to V_{ {\boldsymbol{\sigma}}}$ and contravariant injective actions
 upon tails and edges:
\[
 f^t:\, \bf{T}_{{\boldsymbol{\sigma}}}\to \bf{T}_{\boldsymbol{\tau}}, \quad f^e:\, E_{\boldsymbol{\sigma}}\to E_{\boldsymbol{\tau}}.
\]
 Geometrically, $f$ contracts edges from $E_{{\bf\tau}}\setminus f^e(E_{{\boldsymbol{\sigma}}})$
 and tails from $\bf{T}_{{\boldsymbol{\tau}}}\setminus f^t(\bf{T}_{{\boldsymbol{\sigma}}})$,
 compatibly with its action upon vertices.
 
 \smallskip

Similarly, given a stratum $\overline{D}( {\boldsymbol{\tau}})$ in $B$, its closure is formed from the union of subschemes $D({\boldsymbol{\sigma}})$, such that ${\boldsymbol{\tau}}>{\boldsymbol{\sigma}}$ and, where ${\boldsymbol{\tau}}$ and ${\boldsymbol{\sigma}}$, have the same set of tails. In our case, where the genus of the curve is zero, the condition ${\boldsymbol{\tau}}> {\boldsymbol{\sigma}}$ is uniquely specified by the {\it splitting} data, which can be described as a certain type of Whitehead move.  

\smallskip 

The splitting data is as follows. Choose a vertex $v$ of the set of vertices of $\boldsymbol{\tau}$ and a partition of the set of flags, incident to $v$: $F_{{\boldsymbol{\tau}}}(v)=F'_{{\boldsymbol{\tau}}}(v)\cup F''_{{\boldsymbol{\tau}}}(v)$ such that both subsets are invariant under the involution $j_{\boldsymbol{\tau}}: F_{{\boldsymbol{\tau}}}\to F_{{\boldsymbol{\tau}}}$. To obtain $\boldsymbol{\sigma}$, replace the vertex $v$ by two vertices $v'$ and $v''$ connected by an edge $e$, where the flags verify $F'_{{\boldsymbol{\tau}}}(v')=F'_{{\boldsymbol{\tau}}}(v)\cup \{e'\}$, $F'_{{\boldsymbol{\tau}}}(v'')=F''_{{\boldsymbol{\tau}}}(v)\cup \{e''\}$, where $e',e''$ are the two halves of the edge $e$. The remaining vertices, flags and incidence relations stay the same for ${\boldsymbol{\tau}}$ and $ {\boldsymbol{\sigma}}$. For more details, see ~\cite{Ma99} ch.III $\S$ 2.7, p.90.

 \smallskip

 Finally, the scheme $\M_{0,n}$ is decomposed into pairwise disjoined locally closed strata, indexed by the isomorphism classes of $n$-graphs. Such as depicted in~\cite{Ma99}  ch.III $\S 3$,  the stratification of the scheme is given by trees. If the set $\bf{T}$ is finite, then $\mathcal{T}((\bf{T}))$ is the set of isomorphism classes of trees ${\boldsymbol{\tau}}$, whose external edges are labeled by the elements of $T$. The set of trees is graded by the number of edges:
\[\mathcal{T}((\bf{T}))=\bigcup_{i=0}^{|T|-3}\mathcal{T}_{i}(({\bf T})), \]
where $\mathcal{T}_{i}((\bf{T}))$ is a tree with $i$ edges. The tree $\mathcal{T}_{0}((\bf{T}))$ is the tree with one vertex and the set of flags equals to $\bf{T}$. 

\subsection{Real stratification of the parametrizing space}
\medskip

Compared to the previous exposition, some additional material is necessary to define the real stratification of $\M_{0,n}$. Bridging strata between each other requires the introduction of Whitehead-moves and we recall some tools from real algebraic geometry \cite{ArVa,Gu74,GCR,GCS,San,Vi90}.

\smallskip

In what follows we first focus on the real (resp. imaginary) part of a complex polynomial in $\dpol_n$. This leads to the consideration of harmonic polynomials which can be interpreted as a given class of real algebraic curves in the plane. For simplicity, we consider affine curves (instead of the projective case) and to avoid making any confusions with the notion of drawings, we call a {\it topological model} the level sets of a harmonic polynomial.

\subsection*{Dissipation and perturbation}
{ In the classical literature on real algebraic curves the common terminology to designate a deformation of a given polynomial, and hence of its level curves goes, is the notion of  {\it dissipation} and {\it perturbation}. This comes related to the desingularisation process i.e. one algebraically perturbs an algebraic curve of a given degree in order to desingularise some singular points. In particular, there exists the  Brusotti theorem (Gudkov \cite{Gu74} p.12) which asserts that:

\smallskip 

{\it a singular curve having only nodal points can be perturbed to a curve of the same degree where some singular points are desingularised in an arbitrary prescribed way. }

 More modern developments allow to explicitly perturb the curve in a desired way (see for instance \cite{San} Thm 3.3.)

Regarding the notion of {\it perturbation of an algebraic curve} we refer to works of Gudkov and Viro (see for instance~\cite{Gu74,Vi90}). Let $f$ be an algebraic curve with finitely many singularities.   Let $(f_{\epsilon})$ be a family of algebraic curves defined by polynomials $f_{\epsilon}$ of the same degree as $f=f_0$ and whose coefficients vary continuously with $\epsilon$. For $\epsilon$ close enough to zero the zero-sets of $V(f_{\epsilon})$ are contained in an arbitrarily small neighborhood of $V(f)$ and their topologies coincides with the topology of $V(f)$ except at small neighborhood of the singular points of $f$. The standard way of perturbing for nodal points and singular points of a type $A_1^+$ and $A_t^-$, for $t\geq 3$ and odd in the Viro classification (\cite{Vi90} p. 1098, see also \cite{ArVa}), is by defining $f_{\epsilon}=f+\epsilon g$ where $g$ is a polynomial of the same degree as $f$ and with a finite number of common zeros with $f$ (Lem 4.1 \cite{San}). 

\smallskip 

In order to desingularise any singular points of real algebraic curves, let us invoke the notion of {\it dissipation}. This refers to the change in the topology of a curves in a neighborhood of a singular point.  Let $M_0\subset \mathbb{R}^n$ be an open set and let the function $f:M_0\to \Rl$ be analytic. Denote by $V_U(f)$, where $U\subset {M_0}$, the set $\{x\in U\, |\, f(x)=0\}$. By a singularity of the hypersurface $V_{M_0}(f)$ at the point $x_0$ we mean the class of germs of hypersurfaces which are diffeomorphic to the germ of the hypersurface $V_{M_0}(f)$ at $x_0$.

The function $f$ is included as $f_0$ in a family of analytic functions $f_t:{M_0}\to \Rl$ with $t\in [0,t_0]$. Suppose that this is an analytic family in the sense that the function ${M_0}\times [0,t_0]\to \Rl$ satisfying $(x,t)\mapsto f_t$ is real analytic. 

\smallskip 

Suppose that the hypersurface $V_{M_0}(f_0)$ has an isolated singularity at $x_0$. If there exists a neighborhood $U$ of $x_0$ such that the hypersurfaces $V_{M_0}(f_t)$, $t\in [0,t_0]$ do not have singular points in $U$, then we say that the family of functions $f_t$ {\it dissipates} the singularity $V_{M_0}(f)$ at $x_0$. 

The term of dissipation comes as follows. 
\begin{defn}[Viro~\cite{Vi90},Sec. 3.3, p.1089]
Consider a ball $B\subset \Rl^N$ centered at $x_0$ of small radius, intersecting $V_{M_0}(f)$ at nonsingular points and transversally. Given that $B\subset \Rl^N$ is centered at $x_0$ we require furthermore that:

\begin{itemize}
\item for $t\in [0,t_0]$ the sphere $\partial B$ intersects only transversally $V_{M_0}(f_t)$ at nonsingular points of the hypersurface. 

\item for $t\in (0,t_0]$ the ball $B$ contains no singular points of the hypersurface $V_{M_0}(f_t)$;

\item The pair $(B, V_{\partial B}(f))$ is then homeomorphic to the cone over its boundary $(\partial B, V_{\partial B}(f))$. 
\end{itemize}

Then the family of pairs $(B, V_{B}(f_t))$ with $t\in (0,t_0]$ is called a {\it dissipation} of the germ of the hypersurface $V_{M_0}(f_0)$ at the point $x_0$.
\end{defn}

\begin{rem}When it is clear from context, the ball $B\subset \Rl^N$ centered at $x_0$ is simply denoted $B$ but when it is necessary to specify the radius we write $B_{\epsilon}.$
\end{rem}
\smallskip  

Related techniques of perturbing a polynomial exist namely via (miniversal) unfoldings of singularities. This procedure can be done in the real or complex realm. Given a function $f_0:(\Rl^m,0)\to (\Rl,0)$ with a singular point at 0 (for simplicity), an unfolding $F(x,\lambda)$ of $f_0$ is a germ of a smooth function $F:(\Rl^m\times \Rl^n,0)\to(\Rl,0)$ at the point $(0,0)$ such that $F(x,0)=f_0$. The space $\Rl^n$ of the second argument of $F$ is called the base space of unfolding and $\lambda_1,\cdots,\lambda_n$ are the parameters of the unfolding. 

The deformation is versal if any deformation $G(x,\eta)$ is equivalent to the deformation induced by $F$ (see \cite{ArVa} Chp. 3, p.67). The dimension $n$ of the base space of the versal deformation is the multiplicity of the singularity $f_0$. Such a deformation $F(x,\nu)$ of $f_0$ can be constructed in the following way. Take the quotient ring of the ring of germs of holomorphic functions on $(\Rl^n,0)$ by the ideal generated by the partial derivatives of the function $f$. Then this has the dimension as a real vector space equal to the multiplicity of the singularity of $f_0$.  

}
\subsection*{Desingularisation of double points}
A basic topological operation on the topological model (and hence on its corresponding polynomial) is the desingularization of some of its double points.  The desingularisation of the topological model at a double point $p$ consists in considering a suitable small open neighborhood $U$ of $p$ and substituting $T\cap U$ for two disjoint open curves in such a way that we get a new model with one double point less. There are exactly two ways of desingularising a double point.  This is illustrated below in the figure \ref{F:cross}. 

\begin{figure}[h]
 \begin{center}
 \begin{tikzpicture}[scale=1]
\node (a) at (-2.5,0) {
\begin{tikzpicture}[scale=0.5]
\newcommand{\degree}[0]{4}
\newcommand{\B}[1]{#1*180/\degree}
\newcommand{\R}[1]{#1*180/\degree }
\draw[black,very thick, name=R3R5](\R{3}:1) .. controls(\R{3}:0.3)and (\R{5}:0.3) .. (\R{5}:1) ;
\draw[black,very thick,name=R1R7](\R{1}:1) .. controls(\R{1}:0.3)and (\R{7}:0.3) .. (\R{7}:1) ;
\end{tikzpicture} 
};
\node (b) at (0,0) {
\begin{tikzpicture}[scale=0.5]
\newcommand{\degree}[0]{4}
\newcommand{\B}[1]{#1*180/\degree}
\newcommand{\R}[1]{#1*180/\degree }
\draw[black, very thick, name=R1R5](\R{1}:1) .. controls(\R{1}:0.3)and (\R{5}:0.3) .. (\R{5}:1) ;
\draw[black, very thick,name=R3R7](\R{3}:1) .. controls(\R{3}:0.3)and (\R{7}:0.3) .. (\R{7}:1) ;
\end{tikzpicture} 
};
\node (c) at (2.5,0) {
\begin{tikzpicture}[scale=0.5]
\newcommand{\degree}[0]{4}
\newcommand{\B}[1]{#1*180/\degree}
\newcommand{\R}[1]{#1*180/\degree }
\draw[black, very thick, name=R1R3](\R{1}:1) .. controls(\R{1}:0.3)and (\R{3}:0.3) .. (\R{3}:1) ;
\draw[black, very thick,name=R5R7](\R{5}:1) .. controls(\R{5}:0.3)and (\R{7}:0.3) .. (\R{7}:1) ;
\end{tikzpicture} 
};
\node at (-4,0) {};
 \draw[black,very thick,<->] (a) -- (b);
  \draw[black,very thick,<->]  (b) -- (c);
\end{tikzpicture}
\end{center}
\vspace{-10pt}
\caption{Desingularizing a double point}\label{F:cross}
\end{figure}

In particular, the hidden operation to perform for such an operation is to proceed to a ``ghost operation'' (an operation that does not exist from the point of view of real algebraic curves but only from a topological insight) by inserting a diagonal joining the two branches which have been inserted. 
Those artificial diagonals have been introduced since the topological models are considered modulo topological equivalence and thus one can transform them by global homeomorphisms (see~\cite{San} p.295).

\subsection*{A. Contracting Whitehead moves}
We are interested not only in nodal points but also in singular points ($k$-fold points). To do so,  
we introduced a generalisation of the construction above which we call {\it Whitehead moves}. This allows a more precise description of what is happening in our situation.

Consider a topological model, being an embedded  forest. The graph has an even number of leaves and vertices are of even valency. Suppose that there exists a face (1-connected region in the complement of the embedded graph) of degree $m\geq 2$. Add one vertex (different from the leafs) to each connected component in the boundary of the face and draw a polygon, joining these vertices. Contract this polygon to a point.  
This contracting morphism deforms the first graph by defining a new tree, which is locally star-like. We call this operation a {\it complete contracting Whitehead move}. This enters a  {\it singularisation} process, where a new $m$-fold point arises in the topological model and also on the corresponding perturbed harmonic polynomial. 

\smallskip

This type of operation also holds in the case of a tree where two vertices are connected to each other by an edge (and leads to the nodal singular points). In this situation, the edge can be contracted, leaving only one vertex, exactly as in the contracting morphism discussed above. We call this a {\it partial contracting Whitehead move}.

\smallskip
\begin{figure}[h]
\begin{center}
\begin{tikzpicture}[scale=1]
\node (a) at (-2.5,0) {
\begin{tikzpicture}[scale=1.5]
\newcommand{\degree}[0]{3}
\newcommand{\last}[0]{5}
\newcommand{\B}[1]{#1*180/\degree}
\newcommand{\R}[1]{#1*180/\degree + 90/\degree}
\draw[blue,   name =B0B1] (\B{0}:1) .. controls(\B{0}:-0.2) and (\B{1}:-0.2) .. (\B{1}:1) ;
\draw[blue,  name =B2B3] (\B{2}:1) .. controls(\B{2}:-0.2) and (\B{3}:-0.2) .. (\B{3}:1) ;
\draw[blue,  name=B4B5] (\B{4}:1) .. controls(\B{4}:-0.2) and (\B{5}:-0.2) .. (\B{5}:1) ;
\draw[red, name =R0R1] (\R{0}:1) .. controls(\R{0}:0.1) and (\R{1}:0.1) .. (\R{1}:1) ;
\draw[red, name =R2R3] (\R{2}:1) .. controls(\R{2}:0.1) and (\R{3}:0.1) .. (\R{3}:1) ;
\draw[red, name =R4R5] (\R{4}:1) .. controls(\R{4}:0.1) and (\R{5}:0.1) .. (\R{5}:1) ;
\end{tikzpicture}
};
\node (b) at (2.5,0) {
\begin{tikzpicture}[scale=1.5]
\newcommand{\degree}[0]{3}
\newcommand{\last}[0]{5}
\newcommand{\B}[1]{#1*180/\degree}
\newcommand{\R}[1]{#1*180/\degree + 90/\degree}
\draw[blue,  name=B0B3] (\B{0}:1) .. controls(\B{0}:-0.2) and (\B{3}:-0.2) .. (\B{3}:1) ;
\draw[blue,  name=B1B4] (\B{1}:1) .. controls(\B{1}:-0.2) and (\B{4}:-0.2) .. (\B{4}:1) ;
\draw[blue,  name=B2B5] (\B{2}:1) .. controls(\B{2}:-0.2) and (\B{5}:-0.2) .. (\B{5}:1) ;
\draw[red, name=R0R1] (\R{0}:1) .. controls(\R{0}:0.2) and (\R{1}:0.2) .. (\R{1}:1) ;
\draw[red, name=R2R3] (\R{2}:1) .. controls(\R{2}:0.2) and (\R{3}:0.2) .. (\R{3}:1) ;
\draw[red, name=R4R5] (\R{4}:1) .. controls(\R{4}:0.2) and (\R{5}:0.2) .. (\R{5}:1) ;
\end{tikzpicture}
};
\draw[->,thick] (a) -- (b);
 \end{tikzpicture}
\end{center}
\vspace{-10pt}
\caption{Example of a complete contracting Whitehead move}\label{WH1.pdf}
\end{figure}
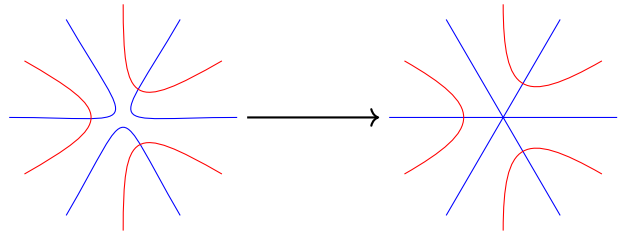

\begin{figure}[h]
\begin{center}
 \begin{tikzpicture}[scale=1]
\node (a) at (-2.5,0) {
\begin{tikzpicture}[scale=1.5]
\newcommand{\degree}[0]{4}
\newcommand{\last}[0]{7}
\newcommand{\B}[1]{#1*180/\degree}
\newcommand{\R}[1]{#1*180/\degree + 90/\degree}
\draw[blue,  name=B0B4] (\B{0}:1) .. controls(\B{0}:0.3) and (\B{4}:0.3) .. (\B{4}:1) ;
\draw[blue,  name=B1B7] (\B{1}:1) .. controls(\B{1}:0.3) and (\B{7}:0.3) .. (\B{7}:1) ;
\draw[blue,  name=B2B6] (\B{2}:1) .. controls(\B{2}:0.3) and (\B{6}:0.3) .. (\B{6}:1) ;
\draw[blue,  name=B3B5] (\B{3}:1) .. controls(\B{3}:0.3) and (\B{5}:0.3) .. (\B{5}:1) ;
\draw[red, name=R0R1] (\R{0}:1) .. controls(\R{0}:0.1) and (\R{1}:0.1) .. (\R{1}:1) ;
\draw[red, name=R2R3] (\R{2}:1) .. controls(\R{2}:0.1) and (\R{3}:0.) .. (\R{3}:1) ;
\draw[red, name=R4R5] (\R{4}:1) .. controls(\R{4}:0.1) and (\R{5}:0.1) .. (\R{5}:1) ;
\draw[red, name=R6R7] (\R{6}:1) .. controls(\R{6}:0.1) and (\R{7}:0.1) .. (\R{7}:1) ;
\end{tikzpicture}
};
\node (b) at (2.5,0) {
\begin{tikzpicture}[scale=1.5]
\newcommand{\degree}[0]{4}
\newcommand{\last}[0]{7}
\newcommand{\B}[1]{#1*180/\degree}
\newcommand{\R}[1]{#1*180/\degree + 90/\degree}
\draw[blue,  name=B0B4] (\B{0}:1) .. controls(\B{0}:0.3) and (\B{4}:0.3) .. (\B{4}:1) ;
\draw[blue,  name=B1B5] (\B{1}:1) .. controls(\B{1}:0.3) and (\B{5}:0.3) .. (\B{5}:1) ;
\draw[blue,  name=B2B6] (\B{2}:1) .. controls(\B{2}:0.3) and (\B{6}:0.3) .. (\B{6}:1) ;
\draw[blue,  name=B3B7] (\B{3}:1) .. controls(\B{3}:0.3) and (\B{7}:0.3) .. (\B{7}:1) ;
\draw[red, name=R0R1] (\R{0}:1) .. controls(\R{0}:0.3) and (\R{1}:0.3) .. (\R{1}:1) ;
\draw[red, name=R2R3] (\R{2}:1) .. controls(\R{2}:0.3) and (\R{3}:0.3) .. (\R{3}:1) ;
\draw[red, name=R4R5] (\R{4}:1) .. controls(\R{4}:0.3) and (\R{5}:0.3) .. (\R{5}:1) ;
\draw[red, name=R6R7] (\R{6}:1) .. controls(\R{6}:0.3) and (\R{7}:0.3) .. (\R{7}:1) ;
\end{tikzpicture}
};
\draw[->,thick] (a) -- (b);
 \end{tikzpicture}
\end{center}
\vspace{-10pt}
\caption{Example of two partial contracting Whitehead moves (done simultaneously)}\label{WH2}
\end{figure}

\subsection*{B. Smoothing Whitehead moves}
The opposite of this operation exists too. It is called a {\it smoothing operation}. Geometrically speaking, a smoothing is applied to the vertex (not a leaf) of a tree having $m\geq2$ edges. This smoothing is 
obtained by ungluing those edges, giving $m$  trees (possibly one-edged) in the boundary of a face of degree $m$. This is called a  {\it complete smoothing Whitehead move}. In this smoothing process, we include also the splitting operation above, i.e. for a given tree with $m\geq3$ edges,  replace one vertex by two vertices connected by an edge, such that the valency of the new vertices remains even, and that the condition for flags holds (where ''flags'' are replaced here by incident half-edges). This is called the {\it partial smoothing Whitehead move}.

\smallskip

\subsection*{Deformation retract Lemma}
The Whitehead moves can be analytically interpreted using the polynomials and the space of their critical values.
 Let $P\in \dpol_{n}$, explicitly: $P(z)=z^{n}+a_{n-1} z^{n-1}+\ldots + a_0$. Let $\tilde{\C}_n$ denote the (affine) space of the critical values $\uv$'s and recall that there is a ramified cover 
$$\pi_w :\dpol_{n} \rightarrow \tilde{\C}_n$$ 
of degree $(n)^{n-1}$~\cite{Myc69}. 

\smallskip 
\begin{thm}\label{polyhedral_complex}
The map $\nu$ that sends a polynomial in $\dpol_{n}$ to its critical values realizes $\dpol_{n}$ as a finite ramified cover of $\tilde{\C}_n$.
\end{thm}
\begin{proof}
The image of $\nu$ contains only unordered tuples of $n-1$ complex numbers different from zero, since a polynomial can have 0 as a critical value if and only if it has multiple roots. Therefore, the image of $\nu$ lies in $\tilde{\C}_n$. To show that $\nu$ is surjective, we use a theorem of R. Thom~\cite{Th63} (1963), stating that given $n-1$ complex critical values, there exists a complex polynomial $P$ of degree $n$ such that $P(r_{i})=v_{i}$ for $1\leq i\leq n-1$, where the $r_{i}$ are the critical points of $P$, and $P(0)=0$.  To find a Tschirnhausen representative polynomial of $\dpol_{n}$ having the same property it suffices to take $P(z-\frac{a_{n-1}}{n-1})$ where, $a_{n-1}$ is the coefficient of $z^{n-1}$ in $P$. By a result of J. Mycielski~\cite{Myc69}, the map $\nu$ is a finite ramified cover, of degree $\frac{n^{n-1}}{n-1}$, see~\cite{Myc69}.
\end{proof}

Let $\mathbb{C}_r^n$ denote the affine space of the critical points $\ur$
and let $p : \mathbb{C}_r\rightarrow \tilde{\C}_n$ denote the natural map given by $P$ ($P(r)=v)$. 
\begin{rem}
Then for a given critical value $v_0$ lying on the real or imaginary line, the smoothing Whitehead move results in the fact that this critical value $v_0$ moves into the interior of one of the quadrants $A,B,C$ or $D$. If one has only a partial smoothing Whitehead move, then the critical value remains  on the half line in the boundary of one of the quadrants, somehow its multiplicity shrinks.  
\end{rem}

We now introduce the deformation retract lemma, which relates to this and is useful for the proof of the fact that the stratification is a Goresky--MacPherson one. 

\begin{lem}[Deformation retract]\label{L:DeR}
Consider a signature $\sigma$. Suppose that in $\sigma$ there exist a face of degree $m$. 
Apply a complete contraction Whitehead move onto the boundary of the face and obtain the signature $\tau$. Then, this operation corresponds to a deformation retract of $\sigma$ onto a new signature: $\tau$. 
\end{lem}
To prove this lemma, we essentially need notions from Morse theory (see \cite{Mi1}) on singular spaces \cite{Laz} and a construction very close to the Milnor fiber. 

Standard Morse theory is concerned with the study of a compact manifold, via the critical points of a Morse function $f$, which is obtained by analysing sublevel sets $F(a)=f^{-1}(\infty,a]$, where $a\in \mathbb{R}$ and how their topology modifies as the point $a$ varies.

\smallskip 

Consider $M$ a real analytic space (of dimension $n$) with isolated singularities. For simplicity we suppose $M$ embedded in some Euclidean space $\mathbb{R}^N$. A Morse function on $M$ is a proper differentiable map $f:M\to \mathbb{R}$ with only non degenerate critical points. We call a critical point $p$ non degenerate iff $p$ is a simple point of $M$ and the Hessian is non singular. One easily check that if $f$ is a Morse function then the set of non-degenerate critical points are isolated. 

One can restrict investigations for $M$ on the local neighbourhood of a given singular point of $M$.  
Suppose $p$ is an isolated singularity of $M$ and $f:M\to \mathbb{R}$ a differentiable function regular at $p$, such that $f(p)=0.$
Then, for a real number $0<\epsilon_{0}<1$ one can define a space with isolated singularities, whose only singularity is $p$ and which is given by the intersection $B_{\epsilon_0}\cap M$, where $B_{\epsilon_0}$ is an open ball of radius $\epsilon_{0}$ centred at $p$.

\smallskip

We recall an important result. 
\begin{prop}[ \cite{Laz}, Prop. 6 . ii]
For $a\in \mathbb{R}$ let $M^a=f^{-1}(-\infty, a]=\{x\in M: f(x)\leq a\}$ be a sublevel set. Let $p$ be a singular point in $M$ such that $f(p)=c$. Moreover, suppose that $a<c<b$ such that the pre-image $f^{-1}([a,b])-\{p\}$ does not contain any singular points of $M$ nor any critical points of $f.$
Then $M^c$ is a deformation retract of $M^b$. 
\end{prop}


\medskip

We are now able to prove the  lemma.
\begin{proof}[Proof of deformation retract lemma \ref{L:DeR}] { 
Let $P_0$ be a complex polynomial in $\dpol_n$.
Consider a given drawing of signature $\tau$, where by hypothesis there exists a (singular) $m$-fold point $z_0$ given by the intersection of $m \geq 2$ curves. Without loss of generality, suppose that this point $z_0$ is an isolated singularity of the real algebraic curve $f_0=\Re P_0(z_0) =0$. 

We construct the germ of the analytic space $M$ using one key ingredient: 
the family of analytic functions $f_t:\mathbb{R}^2\to \Rl$ with $t\in [0,t_0]$ where we have a real analytic germ 
given by $F(x,t):M\times [0,t_0]\to \Rl$ such that it dissipates the isolated singularity given by $f_0$ at $z_0=0$.
In particular by definition of the dissipation, we consider this locally within a ball $B$ of small radius $\epsilon$ centered at $z_0$. 
$M$ can be topologically interpreted as the cone over the pair $(\partial B,V_{\partial B}(f_0))$. For every $t\in [0,t_0]$ the topological model of $f_t=0$ 
corresponds to the level curves on the surface $M$. In particular, for $t=0$, there exists a singular point $p$ on $M$; for $t\neq 0$, one has desingularised the $m$-fold point by a complete smoothing Whitehead operation.

}
So, to summarise we have constructed a (real) germ $M$ with an isolated singularity at 0 of multiplicity $m$.

It remains to construct the Morse function $d:M\to \mathbb{R}$. This follows almost immediately if we take the square of the function ``distance to 0''  i.e.$d(x_1,\cdots,x_N)=x_1^2+\cdots + x_N$, where $N$ is the dimension in which the ball $B$ is defined.
The existence of the critical point  in $M$ is given at $0$ and so we can put $d(0)=0$. 

Solutions of its gradient vector field $\nabla d$ are the straight rays that emanate from the origin. We adapt this vector field to $M$. At each point $x\in M\setminus \{0\}$ the gradient vector is obtained by projecting $\nabla d_M(x)$ to $T_x(M\setminus \{0\})$ (the tangent space to $M\setminus \{0\}$). 

This means that a point $x \in M\setminus \{0\}$ is a critical point of $d_M$ iff $M\setminus \{0\}$ is tangent at $x$ to the sphere passing through $x$ and centered at 0. 
One has that $d_M$ has at most a finite number of critical values corresponding to points in $M\setminus \{0\}$, since it is the restriction of an analytic function on $B$. Hence $M\setminus \{0\}$ meets transversally all sufficiently small spheres around the origin in $\Rl^N$. The gradient vector field of $d_M$ is now everywhere transversal to the spheres around 0, and it can be assumed to be integrable. Hence it defines a one-parameter family of local diffeomorphisms of $M\setminus \{0\}$  taking each link (i.e. $\partial B_{\epsilon}\cap M$) into ``smaller'' links $\partial B_{t-\nu}\cap M$, where $\nu\in(0,\epsilon]$.

One needs to choose real numbers $a,b$ such that for $a\in \mathbb{R}$ (or $b\in \mathbb{R}$) the level subset $M^a=d^{-1}(-\infty, a]=\{x\in M:\, d(x)\leq a\}$ (resp. $M^b$) are consistently defined. So, we choose $a,b$ such that $a<c<b$ (for us $c$ here is given by $c=0$). Moreover, $d$ is constructed such that there are no other singular points in $d^{-1}([a,b])-\{p\}$ nor any other critical points of $d$. 
So, applying the above proposition, $M^c$ (where $c=0$) is a deformation retract of $M^b.$

\end{proof}

\smallskip

\subsection{Some properties of the stratification}

\begin{lem}~\label{L:smo}Suppose $P$ is a polynomial with signature $\tau$ and a critical point of multiplicity $m$ at $z_0$, and let $\sigma$ be the signature obtained by smoothing $\tau$. Then $A_\tau \subset \overline{A}_\sigma$ and there exists a neighborhood $U$ of $P$ in $A_\sigma$ such that $U\cong V \times (\mathbb{D})^m$, where $V$ is a neighborhood of $P\in A_\tau$ and the polydisk $(\mathbb{D})^m$ corresponds to a canonical local perturbation of $(z - z_0)^m$.

\end{lem}

The proof relies on the perturbation of polynomials and here it applies as well for the complex case as for the real case i.e. for critical points on the real (resp. imaginary) part of the polynomial and for the multiple roots (so on the discriminant variety). 
We start with the easier case which is the latter one.

\begin{proof}
Consider a polynomial $P_0\in A_{\tau}$. Let $V$ be an open neighbourhood of $P_0$ in the stratum $A_\tau$. We have that degree $n$ polynomials in $V$ are of the following shape $(z-z_0)^mg(z)$, where $g$ is a univariate complex polynomial of degree $n-m$. One needs to desingularise the multiplicity $m$ root through a complete smoothening Whitehead move in order to obtain a deformed polynomial in $U$.
 Analytically speaking, there are two operations which occur at the same time:
 \begin{enumerate} 
 \item one dissipates the singular point $z_0$ by perturbing the coefficients of $P_0$ in such a way that we define a family $(P_{\epsilon})$, $0\leq \epsilon\leq 1$ where $P_0\in V$ and $P_1\in U$; in particular, here it means that the topologies of the perturbed level sets of $\Re P_0=0$ and $\Im P_0=0$ coincide with the ones of $\Re P_1=0$ and $\Im P_1=0$ except at a small neighbourhood of the singular point $z_0$. 
 \item one desingularises the singular point $z_0$ using locally a miniversal deformation at the singular germ describing the singularity at $z_0$. \end{enumerate}
The combination of  those two steps together allows the definition of a dissipation of the singular point $z_0$.

Consider a small neighborhood of the singular point $z_0$. Locally, it can be described as a germ with a singularity at the origin of the following type: $(z_0^m,0)$. Regarding the desingularisation of $(z-z_0)^m$, we cannot explicitly write down a universal family in terms of the coefficients $a_i$ of the polynomial $P_0$. However, we do know that there {\it exists} such a  local universal family and that it is biholomorphically equivalent to the one obtained by completing
the polynomial $(z-z_0)^m$ into the generic polynomial in $(z-z_0)$ of degree $m$ near $z_0$ (for instance, see~\cite{V90}, chapter 2, paragraph 3 for more details). In other words, let $\ueps=(\eps_1,\ldots, \eps_m)$ be complex numbers with $\vert \eps_i \vert <\eps<\!< 1$ 
for all $i$ (with some $\eps >0$) and let 
$$p_\eps(z) = (z-z_0)^m+\eps_1(z-z_0)^{m-1}+\ldots+\eps_m \ .$$
Then, there is a biholomorphic map between the set of the $\ueps$'s (i.e. the polydisk $\mathbb{D}_\eps^m$), or equivalently the 
family $p_\eps$, and a universal family $P_\eps(z)$ with $P_0=P$.
Therefore, given a small neighbourhood $U$ in $A_{\sigma}$ one can identify this to polynomials in $V$ for which an $m$-parameter unfolding has been given at the singular point $z_0$ (and the latter is provided by the parameter space $\mathbb{D}_\eps^m$).  

\smallskip 

Concerning a real version of this (i.e. the critical point is different from zero and has a real or imaginary critical value) consider a harmonic polynomial, say $\Re P_0=0$ (resp. $\Im P_0=0$). Suppose that there exists a multiplicity $m$ critical point on this real algebraic curve. Consider a small neighbourhood of this point. Similarly to the previous approach, we dissipate this singular point by defining a path in the space $\dpol_n$ where coefficients of $P_0$ (and thus of $\Re P_0$) are perturbed in such a way that  the critical point $z_0$ gets desingularised (and there are no other modifications concerning the set of singular points of $\Re P_0, \Im P_0=0$  and thus $P_0$). 

However, similarly to the previous case, one needs to consider the germ of the singularity at the point $z_0$ of $\Re P_0=0$ and deform it in order to define the dissipation. Suppose that locally the singularity at $z_0$ is given by the real germ $f:(\Rl^2,0)\to (\Rl,0)$. Then, one can unfold the singularity using the family of mappings $F:\Rl^2\times\Rl^m\to \Rl$, where $\Rl^m$ is the $m$-parameter unfolding of $f$ (the dimension of the parameter space is given by the multiplicity of the critical point which is $m$, here). 

Again, we cannot explicitly write down a universal family in terms of the coefficients $a_i$ of the polynomial $P_0$ (and thus of coefficients of $\Re P_0$). But we do know that there exists such a local universal family and that it is homeomorphically equivalent to the one obtained by completing the polynomial $f$ defining the germ of the singularity at $z_0$.  So, by doing this we have dissipated the singular point. This construction allows to state that in a small neighbourhood $U$ in $A_{\sigma}$ one can identify it as $V\times \mathbb{D}_\eps^m$ where $V$ is a small neighborhood in $A_{\tau} $ and $\mathbb{D}_\eps^m$ is the parameter space allowing an unfolding of the singular point. 

\end{proof}

\begin{lem}\label{L:6}
Strata in the topological closure $\overline{A_{\sigma}}$ are indexed by signatures obtained from a contracting Whitehead  move on $\sigma$. 
\end{lem}

\begin{proof}
Let us suppose first that we apply on a generic signature $\sigma$ a complete contracting Whitehead.
Let $P$ be a polynomial in $A_{\sigma}$ and $\tilde{P}$ a polynomial in $A_{\tilde{\sigma}}$, where $\sigma\prec \tilde{\sigma}$. 
We first use the Whitehead move of first type. 
A contracting Whitehead move on $\sigma$ corresponds to a path of the critical values of $P$ in  the space of critical values. We will show that by using a contracting Whitehead move on $\sigma$ this corresponds to defining a convergent sequence of critical values in $ \mathbb{C}_w^{n-1}$. 
In $\tilde{\sigma}$ the intersection point of a set of $m$ curves of the same color lies on a critical point $c_i$ of $\tilde{P}$. Suppose that there exist $I$ (where $|I|<d$) such intersections. Therefore, 
for $i\in I$ we have $\tilde{P}'(r_i)=0$ and $\Re(\tilde{P}(r_i))=0$ (resp $\Im(\tilde{P}(r_i))=0$) and the critical value $\tilde{P}(r_i)=v_i  \in \imath\mathbb{R}$ (resp. $\tilde{P}(r_i)=v_i  \in \mathbb{R}$).  
So, this indicates a sequence of critical values converging to  $(v_1,....,v_{n-1})$, where the subset $v_i, i \in \{1,...,n-1\}$  lies on the imaginary (resp. real) axes.
Hence it indicates a topological closure. 

\smallskip

Consider the case of a partial contracting Whitehead move.
In this case, the initial signature $\sigma$ has a set of critical values lying on the real or imaginary axis and the partial contracting Whitehead operation merges a subset of those critical values together. A partial contracting Whitehead move corresponds to a converging sequence of critical values in the space $\mathbb{C}_v^{n-1}$ and hence it indicates a topological closure of the stratum $A_{\sigma}$. 
\end{proof}

\subsection{$\dpol_{d}$ as a covering of a non-compact stratified space}\label{S:crit}
\subsubsection{Critical values}

Consider the space $V_{n}=(\mathbb{C}^{n-1}\setminus 0)/S_{n-1}$, where $S_{n-1}$ is the group of permutations. If $X$ denotes an equivalence class
of points in $\mathbb{C}^{n-1}$, we can associate a unique $\sigma$-sequence $(a,b,c,d,e,f,g,h)$ of positive integers to $X$ enumerating the number of points 
in $X$ in the quadrants $A,B,C,D$ and on the semi-axes. The set of points $X$ in $V_{n}$ having a given $\sigma$-sequence $(a,b,c,d,e,f,g,h)$  forms a polygonal cell in $V_{d}$ isomorphic to 

\begin{equation}\label{E:poly}
 A^{a} /S_{a} \times B^{b}/ S_{b}\times C^{c}/ S_{c} \times D^{d}/ S_{d} \times (\mathbb{R}^{+})^{e} /S_{e} \times (\mathbb{R}^{-})/ S_{f}\times  (\imath\mathbb{R}^{+})^{g}/ S_{g} \times (\imath\mathbb{R}^{-})^{h}/ S_{h}.
 \end{equation}
The real dimension of this cell is equal to $2(a+b+c+d)+ (e+f+g+h)$. The cells are disjoint and thus form a stratification of $V_{n}$.   
\begin{defn}
 A subset $V$ inside $\mathbb{C}^{n-1}/S_{n-1}$ for $n>2$ is said to be a $\mathbf {non-compact\ \ stratification}$ if it is equipped with a stratification by a finite number of open cells of varying dimensions having the following properties:
 \begin{itemize}
 \item the relative closure of a $k$-dimensional cell of $V$ is a union of cells in the stratification. 
 \item the relative closure of a $k$-dimensional cell of $V$ is a ``semi-closed'' polytope, i.e. the union of the interior of a closed polytope in $\mathbb{C}^{n-1}/S_{n-1}$ with a subset of its faces. 
 \end{itemize} 

\end{defn}
\begin{lem}
Let $n>2$, and let $\mathcal{V}_{n}$ denote the space $V_{n}$ equipped with the stratification by $\sigma$-sequences. Then $\mathcal{V}_{n}$ is a non-compact stratification. 
\end{lem}
\begin{proof}
The closure of the region of $\mathcal{V}_{n}$ described by~\ref{E:poly} is given by $$\overline{A}^{a} /S_{a} \times \overline{ B}^{b}/ S_{b}\times \overline{ C}^{c}/ S_{c} \times \overline{D}^{d}/ S_{d} \times (\mathbb{R}^{+})^{e} /S_{e} \times (\mathbb{R}^{-})/ S_{f}\times  (\imath\mathbb{R}^{+})^{g}/ S_{g} \times (\imath\mathbb{R}^{-})^{h}/ S_{h},$$ where $\overline{A}$ denotes the closure in $V_1=\mathbb{C}\setminus 0$ of the quadrant $A$, namely the union of $A$ with $\mathbb{R}^{+}$ and $\imath\mathbb{R}^{+}$, and similarly for  $\overline{B}, \overline{C}, \overline{D}$. 
The direct product of semi-closed polytopes is again a semi-closed polytope, as is the quotient of a semi-closed polytope by a sub-group of its symmetry group. 
 \end{proof}

\subsection{The cases $n=2,3,4$}
The exact nature of the ramified cover $\dpol_{n} \rightarrow \mathcal{V}_{n}$ is complicated and interesting, especially in terms of describing the ramification using the signatures.  In this section, we work out full details in the small dimensional cases, and for generic strata.

\begin{figure}[h]
\begin{center}
 \begin{tikzpicture}[scale=0.8]
\node (a) at (-2,2) {
 \begin{tikzpicture}[scale=0.7]
\newcommand{\degree}[0]{8}
\newcommand{\last}[0]{7}
\newcommand{\B}[1]{#1*360/\degree}
\draw[black,thick] (0,0) circle (1) ;
\draw[blue, thick, name=B0B6](\B{0}:1) .. controls(\B{0}:0.3)and (\B{6}:0.3) .. (\B{6}:1) ;
\draw[blue, thick, name=B2B4](\B{2}:1) .. controls(\B{2}:0.3)and (\B{4}:0.3) .. (\B{4}:1) ;
\draw[red,thick, name=B1B3](\B{1}:1) .. controls(\B{1}:0.3)and (\B{3}:0.3) .. (\B{3}:1) ;
\draw[red,thick, name=B5B7](\B{5}:1) .. controls(\B{5}:0.3)and (\B{7}:0.3) .. (\B{7}:1) ;
\end{tikzpicture} 
};
\node (b) at (2,2) {
\begin{tikzpicture}[scale=0.7]
\newcommand{\degree}[0]{8}
\newcommand{\last}[0]{7}
\newcommand{\B}[1]{#1*360/\degree}
\draw[black,thick] (0,0) circle (1) ;
\draw[blue, thick, name=B0B6](\B{0}:1) .. controls(\B{0}:0.3)and (\B{6}:0.3) .. (\B{6}:1) ;
\draw[blue, thick, name=B2B4](\B{2}:1) .. controls(\B{2}:0.0)and (\B{4}:0.3) .. (\B{4}:1) ;
\draw[red, thick,  name=B1B7](\B{1}:1) .. controls(\B{1}:0.3)and (\B{7}:0.3) .. (\B{7}:1) ;
\draw[red, thick, ,name=B3B5](\B{3}:1) .. controls(\B{3}:0.3)and (\B{5}:0.3) .. (\B{5}:1) ;
\end{tikzpicture} 
};
\node (c) at (2,-2) {
\begin{tikzpicture}[scale=0.7]
\newcommand{\degree}[0]{8}
\newcommand{\last}[0]{7}
\newcommand{\B}[1]{#1*360/\degree}
\draw[black,thick] (0,0) circle (1) ;
\draw[blue, thick, name=B0B2](\B{0}:1) .. controls(\B{0}:0.3)and (\B{2}:0.3) .. (\B{2}:1) ;
\draw[blue, thick, name=B4B6](\B{4}:1) .. controls(\B{4}:0.3)and (\B{6}:0.3) .. (\B{6}:1) ;
\draw[red, thick, name=B3B5](\B{3}:1) .. controls(\B{3}:0.3)and (\B{5}:0.3) .. (\B{5}:1) ;
\draw[red, thick,name=B1B7](\B{1}:1) .. controls(\B{1}:0.3)and (\B{7}:0.3) .. (\B{7}:1) ;
\end{tikzpicture} 
};
\node (d) at (-2,-2) {
\begin{tikzpicture}[scale=0.7]
\newcommand{\degree}[0]{8}
\newcommand{\last}[0]{7}
\newcommand{\B}[1]{#1*360/\degree}
\draw[black,thick] (0,0) circle (1) ;
\draw[blue, thick, name=B0B2](\B{0}:1) .. controls(\B{0}:0.3)and (\B{2}:0.3) .. (\B{2}:1) ;
\draw[blue, thick, name=B4B6](\B{4}:1) .. controls(\B{4}:0.3)and (\B{6}:0.3) .. (\B{6}:1) ;
\draw[red, thick, name=B5B7](\B{5}:1) .. controls(\B{5}:0.3)and (\B{7}:0.3) .. (\B{7}:1) ;
\draw[red, thick, name=B1B3](\B{1}:1) .. controls(\B{1}:0.3)and (\B{3}:0.3) .. (\B{3}:1) ;
\end{tikzpicture} 
};
 \draw[->, black,thick] (a) -- (b);
 \draw[->,black,thick] (b) -- (c);
 \draw[->,black,thick] (c) -- (d);
 \draw[->,black,thick] (d) -- (a);
 \draw (-2,2) node {\Large $1$};     
  \draw (2,2) node {\Large $2$}; 
   \draw (2,-2) node {\Large $3$};     
  \draw (-2,-2) node {\Large $4$}; 
\end{tikzpicture} 
\end{center}
\vspace{-8pt}
\caption{Relations between strata for $n=2$}\label{F:d=2}
\end{figure}
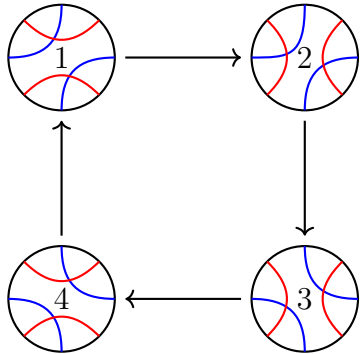

\vskip.1cm
Let $n=2$. The spaces $\dpol_{2}$ and $\mathcal{V}_{2}$ are one-dimensional. The space $\mathcal{V}_{2}$ is $\mathbb{C}\setminus 0$ equipped with the stratification given by the four quadrants $A,B,C,D$ and the four semi-axes. The only Tschirnhausen polynomial of degree 2 having given critical value $v$ is $z^2 +v$, therefore the covering map $\nu$ is unramified of degree 1, an isomorphism. The four signatures corresponding to the strata of real dimensional 2 and the four corresponding to the one dimensional strata are illustrated in figure~\ref{F:d=2}.
  
\smallskip
 
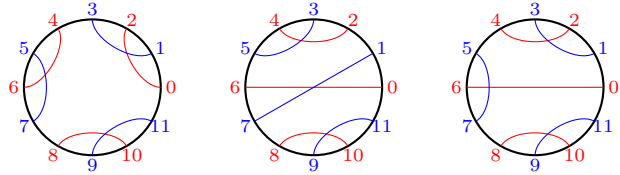
\begin{figure}[h]
\begin{center}
\begin{tikzpicture}[scale=0.9]
\newcommand{\degree}[0]{6}
\newcommand{\B}[1]{#1*180/\degree}
\newcommand{\R}[1]{#1*180/\degree}
\newcommand{\RDOT}[0]{[red](i-1) circle (0.05)}
\newcommand{\BDOT}[0]{[blue](i-1) circle (0.05)}
\draw[black,thick] (0,0) circle (1) ;
\foreach \k in {0,2,4,6,8,10} {
  \draw[red] (\R{\k}:1.15) node {\scriptsize\k};} ;
 \foreach \k in {1,3,5,7,9,11} { 
\draw[ blue] (\B{\k}:1.15) node {\scriptsize\k} ;};
\draw[red, name=B0B2](\R{0}:1) .. controls(\R{0}:0.7)and (\R{2}:0.7) .. (\R{2}:1) ;
\draw[red,name=B4B6](\R{4}:1) .. controls(\R{4}:0.7)and (\R{6}:0.7) .. (\R{6}:1) ;
\draw[red,name=B8B10](\R{8}:1) .. controls(\R{8}:0.7)and (\R{10}:0.7) .. (\R{10}:1) ;
\draw[blue,  name=R1R3](\B{1}:1) .. controls(\B{1}:0.7)and (\B{3}:0.7) .. (\B{3}:1) ;
\draw[blue, name=R5R7](\B{5}:1) .. controls(\B{5}:0.7)and (\B{7}:0.7) .. (\B{7}:1) ;
\draw[blue, name=R9R11](\B{9}:1) .. controls(\B{9}:0.7)and (\B{11}:0.7) .. (\B{11}:1) ;
\end{tikzpicture}
\quad
\begin{tikzpicture}[scale=0.9]
\newcommand{\degree}[0]{6}
\newcommand{\B}[1]{#1*180/\degree}
\newcommand{\R}[1]{#1*180/\degree }
\newcommand{\RDOT}[0]{[red](i-1) circle (0.05)}
\newcommand{\BDOT}[0]{[blue](i-1) circle (0.05)}
\draw[black,thick] (0,0) circle (1) ;
\foreach \k in {0,2,4,6,8,10} {
  \draw[red] (\R{\k}:1.15) node { \scriptsize\k};} ;
  \foreach \k in {1,3,5,7,9,11} {
 \draw[ blue] (\B{\k}:1.15) node {\scriptsize\k} ;
}
\draw[red, name=B0B6](\R{0}:1) .. controls(\R{0}:0.7)and (\R{6}:0.7) .. (\R{6}:1) ;
\draw[red, name=B2B4](\R{2}:1) .. controls(\R{2}:0.7)and (\R{4}:0.7) .. (\R{4}:1) ;
\draw[red, name=B8B10](\R{8}:1) .. controls(\R{8}:0.7)and (\R{10}:0.7) .. (\R{10}:1) ;
\draw[blue, name=R1R7](\B{1}:1) .. controls(\B{1}:0.7)and (\B{7}:0.7) .. (\B{7}:1) ;
\draw[blue, name=R3R5](\B{3}:1) .. controls(\B{3}:0.7)and (\B{5}:0.7) .. (\B{5}:1) ;
\draw[blue, name=R9R11](\B{9}:1) .. controls(\B{9}:0.7)and (\B{11}:0.7) .. (\B{11}:1) ;
\end{tikzpicture}
\quad
\begin{tikzpicture}[scale=0.9]
\newcommand{\degree}[0]{6}
\newcommand{\B}[1]{#1*180/\degree}
\newcommand{\R}[1]{#1*180/\degree }
\newcommand{\RDOT}[0]{[red](i-1) circle (0.05)}
\newcommand{\BDOT}[0]{[blue](i-1) circle (0.05)}
\draw[black,thick] (0,0) circle (1) ;
\foreach \k in {0,2,4,6,8,10} {
  \draw[red] (\R{\k}:1.15) node { \scriptsize\k};} ;
  \foreach \k in {1,3,5,7,9,11} {
 \draw[ blue] (\B{\k}:1.15) node {\scriptsize\k} ;
}
\draw[red, name=B0B6](\R{0}:1) .. controls(\R{0}:0.7)and (\R{6}:0.7) .. (\R{6}:1) ;
\draw[red, name=B2B4](\R{2}:1) .. controls(\R{2}:0.7)and (\R{4}:0.7) .. (\R{4}:1) ;
\draw[red, name=B8B10](\R{8}:1) .. controls(\R{8}:0.7)and (\R{10}:0.7) .. (\R{10}:1) ;
\draw[blue, name=R1R3](\B{1}:1) .. controls(\B{1}:0.7)and (\B{3}:0.7) .. (\B{3}:1) ;
\draw[blue, name=R5R7](\B{5}:1) .. controls(\B{5}:0.7)and (\B{7}:0.7) .. (\B{7}:1) ;
\draw[blue, name=R9R11](\B{9}:1) .. controls(\B{9}:0.7)and (\B{11}:0.7) .. (\B{11}:1) ;
\end{tikzpicture}
\vspace{-10pt}
\caption{Diagrams $d=3$}\label{F:diagrams3}
\end{center}
\end{figure}

 Let $n=3$. In this case the covering map $\nu$ is of degree 3: explicitly, if $P(z)=z^3+az+b$ has critical values $v_1$ and $v_2$ then so do the polynomials $P(\zeta z)$ and  $P(\zeta^2 z)$ where $\zeta^3 =1$. The 10 open strata of real dimension 4 in $\mathcal{V}_3$ are given by: 
\[A\times A/S_{2},\quad B\times B/S_{2},\quad C\times C/S_{2},\quad D\times D/S_{2}, \]
\[A\times B,\quad A\times C,\quad A\times D, \quad B\times C,\quad  B\times D, \quad C\times D.\]
Ramification occurs only above the first four; in fact exactly when the two critical values are equal, i.e. $P(z)=z^3+b$. Thus above each of the first four cells, there is only one stratum, corresponding to the four rotations of the left most signature (see figure~\ref{F:diagrams3}).

\smallskip  

 In contrast, there are three distinct signatures above each of the remaining 6 strata. The six signatures in the middle of the figure~\ref{F:diagrams3} form two orbits under the $\frac{2\pi}{3}$ rotation which lie above the cells, $ A\times C$ and $ B\times D$,
whereas the twelve signatures on the right form four orbits lying over the strata $A\times B; B\times C; C\times D; D\times A$. 

\smallskip

Let $n=4$. The degree of the covering map $\nu$ is 16. There are 20 generic strata in $\mathcal{V}_4$. Four generic strata correspond to taking three critical values in three different quadrants. There is no ramification above these strata; each of these have 16 distinct strata in the preimage of $\nu$, corresponding to four rotations each of the fifth, sixth, eighth and eleventh signatures in the figure below.   

Four more cells of $\mathcal{V}_4$ correspond to taking three critical values in the same quadrant. Only one stratum of $\dpol_{4}$ lies above each of these cells, namely, the first signature in the figure below, with ramification of order sixteen. 
The remaining twelve cells correspond to two critical points in one quadrant and the third in a different quadrant.
When the quadrants are adjacent, only six cells lie above the corresponding regions of  $\mathcal{V}_4$. For example over the region $A,A,B$ lie the two rotations of the second signature in the figure below with ramification of order 2, and the four rotations of the tenth signature each of ramification of order 3.  
The situation is analogous when the quadrants are opposed. For example over the region $A,A,C$, there are the two rotations of the third figure (below) each with ramification of order 2, and the four rotations of the ninth signature each of ramification order 3.   
 \begin{itemize}
 \item All 1 classes of size 4:
 
 \vspace{10pt}
 
 \begin{tikzpicture}[scale=1]
\newcommand{\degree}[0]{16}
\newcommand{\last}[0]{15}
\newcommand{\B}[1]{#1*360/\degree}
\draw (0,0) circle (1) ;
\draw[blue, name=B0B2] (\B{0}:1) .. controls(\B{0}:0.7) and (\B{2}:0.7) .. (\B{2}:1) ;
\draw[blue, name=B4B6] (\B{4}:1) .. controls(\B{4}:0.7) and (\B{6}:0.7) .. (\B{6}:1) ;
\draw[blue, name=B8B10] (\B{8}:1) .. controls(\B{8}:0.7) and (\B{10}:0.7) .. (\B{10}:1) ;
\draw[blue, name=B12B14] (\B{12}:1) .. controls(\B{12}:0.7) and (\B{14}:0.7) .. (\B{14}:1) ;
\draw[red, name=B1B3] (\B{1}:1) .. controls(\B{1}:0.7) and (\B{3}:0.7) .. (\B{3}:1) ;
\draw[red, name=B5B7](\B{5}:1) .. controls(\B{5}:0.7) and (\B{7}:0.7) .. (\B{7}:1) ;
\draw[red, name=B9B11] (\B{9}:1) .. controls(\B{9}:0.7) and (\B{11}:0.7) .. (\B{11}:1) ;
\end{tikzpicture}

\vspace{5pt}
\item All 3 classes of size 8:

\vspace{10pt}

\begin{tikzpicture}[scale=1]
\newcommand{\degree}[0]{16}
\newcommand{\last}[0]{15}
\newcommand{\B}[1]{#1*360/\degree}
\draw (0,0) circle (1) ;
\draw[blue, name=B0B2] (\B{0}:1) .. controls(\B{0}:0.7) and (\B{2}:0.7) .. (\B{2}:1) ;
\draw[blue, name=B4B6] (\B{4}:1) .. controls(\B{4}:0.7) and (\B{6}:0.7) .. (\B{6}:1) ;
\draw[blue, name=B8B10] (\B{8}:1) .. controls(\B{8}:0.7) and (\B{10}:0.7) .. (\B{10}:1) ;
\draw[blue, name=B12B14] (\B{12}:1) .. controls(\B{12}:0.7) and (\B{14}:0.7) .. (\B{14}:1) ;
\draw[red, name=B1B7] (\B{1}:1) .. controls(\B{1}:0.7) and (\B{7}:0.7) .. (\B{7}:1) ;
\draw[red, name=B3B5] (\B{3}:1) .. controls(\B{3}:0.7) and (\B{5}:0.7) .. (\B{5}:1) ;
\draw[red, name=B9B15] (\B{9}:1) .. controls(\B{9}:0.7) and (\B{15}:0.7) .. (\B{15}:1) ;
\draw[red, name=B11B13] (\B{11}:1) .. controls(\B{11}:0.7) and (\B{13}:0.7) .. (\B{13}:1) ;
\end{tikzpicture}
\quad
\begin{tikzpicture}[scale=1]
\newcommand{\degree}[0]{16}
\newcommand{\last}[0]{15}
\newcommand{\B}[1]{#1*360/\degree}
\draw (0,0) circle (1) ;
\draw[blue, name=B0B2] (\B{0}:1) .. controls(\B{0}:0.7) and (\B{2}:0.7) .. (\B{2}:1) ;
\draw[blue, name=B4B14] (\B{4}:1) .. controls(\B{4}:0.7) and (\B{14}:0.7) .. (\B{14}:1) ;
\draw[blue,  name=B6B12] (\B{6}:1) .. controls(\B{6}:0.7) and (\B{12}:0.7) .. (\B{12}:1) ;
\draw[blue, name=B8B10] (\B{8}:1) .. controls(\B{8}:0.7) and (\B{10}:0.7) .. (\B{10}:1) ;
\draw[red, name=B1B3] (\B{1}:1) .. controls(\B{1}:0.7) and (\B{3}:0.7) .. (\B{3}:1) ;
\draw[red, name=B5B15] (\B{5}:1) .. controls(\B{5}:0.7) and (\B{15}:0.7) .. (\B{15}:1) ;
\draw[red, name=B7B13] (\B{7}:1) .. controls(\B{7}:0.7) and (\B{13}:0.7) .. (\B{13}:1) ;
\draw[red, name=B9B11] (\B{9}:1) .. controls(\B{9}:0.7) and (\B{11}:0.7) .. (\B{11}:1) ;
\end{tikzpicture}
\quad
\begin{tikzpicture}[scale=1]
\newcommand{\degree}[0]{16}
\newcommand{\last}[0]{15}
\newcommand{\B}[1]{#1*360/\degree}
\draw (0,0) circle (1) ;
\draw[blue, name=B0B2] (\B{0}:1) .. controls(\B{0}:0.7) and (\B{2}:0.7) .. (\B{2}:1) ;
\draw[blue, name=B4B6] (\B{4}:1) .. controls(\B{4}:0.7) and (\B{6}:0.7) .. (\B{6}:1) ;
\draw[blue, name=B8B10] (\B{8}:1) .. controls(\B{8}:0.7) and (\B{10}:0.7) .. (\B{10}:1) ;
\draw[blue, name=B12B14] (\B{12}:1) .. controls(\B{12}:0.7) and (\B{14}:0.7) .. (\B{14}:1) ;
\draw[red, name=B1B3] (\B{1}:1) .. controls(\B{1}:0.7) and (\B{3}:0.7) .. (\B{3}:1) ;
\draw[red, name=B5B15] (\B{5}:1) .. controls(\B{5}:0.7) and (\B{15}:0.7) .. (\B{15}:1) ;
\draw[red, name=B7B13] (\B{7}:1) .. controls(\B{7}:0.7) and (\B{13}:0.7) .. (\B{13}:1) ;
\draw[red, name=B9B11] (\B{9}:1) .. controls(\B{9}:0.7) and (\B{11}:0.7) .. (\B{11}:1) ;
\end{tikzpicture}\\

\item All 7 classes of size 16:

\vspace{10pt}
\begin{tikzpicture}[scale=1]
\newcommand{\degree}[0]{16}
\newcommand{\last}[0]{15}
\newcommand{\B}[1]{#1*360/\degree}
\draw (0,0) circle (1) ;
\draw[blue, name=B0B2] (\B{0}:1) .. controls(\B{0}:0.7) and (\B{2}:0.7) .. (\B{2}:1) ;
\draw[blue, name=B4B6] (\B{4}:1) .. controls(\B{4}:0.7) and (\B{6}:0.7) .. (\B{6}:1) ;
\draw[blue, name=B8B14] (\B{8}:1) .. controls(\B{8}:0.7) and (\B{14}:0.7) .. (\B{14}:1) ;
\draw[blue, name=B10B12] (\B{10}:1) .. controls(\B{10}:0.7) and (\B{12}:0.7) .. (\B{12}:1) ;
\draw[red, name=B1B3] (\B{1}:1) .. controls(\B{1}:0.7) and (\B{3}:0.7) .. (\B{3}:1) ;
\draw[red, name=B5B15] (\B{5}:1) .. controls(\B{5}:0.7) and (\B{15}:0.7) .. (\B{15}:1) ;
\draw[red, name=B7B13] (\B{7}:1) .. controls(\B{7}:0.7) and (\B{13}:0.7) .. (\B{13}:1) ;
\draw[red, name=B9B11] (\B{9}:1) .. controls(\B{9}:0.7) and (\B{11}:0.7) .. (\B{11}:1) ;
\end{tikzpicture}
\quad
\begin{tikzpicture}[scale=1]
\newcommand{\degree}[0]{16}
\newcommand{\last}[0]{15}
\newcommand{\B}[1]{#1*360/\degree}
\draw (0,0) circle (1) ;
\draw[blue, name=B0B2] (\B{0}:1) .. controls(\B{0}:0.7) and (\B{2}:0.7) .. (\B{2}:1) ;
\draw[blue, name=B4B6] (\B{4}:1) .. controls(\B{4}:0.7) and (\B{6}:0.7) .. (\B{6}:1) ;
\draw[blue, name=B8B14] (\B{8}:1) .. controls(\B{8}:0.7) and (\B{14}:0.7) .. (\B{14}:1) ;
\draw[blue, name=B10B12] (\B{10}:1) .. controls(\B{10}:0.7) and (\B{12}:0.7) .. (\B{12}:1) ;
\draw[red, name=B3B9] (\B{3}:1) .. controls(\B{3}:0.7) and (\B{9}:0.7) .. (\B{9}:1) ;
\draw[red, name=B5B7] (\B{5}:1) .. controls(\B{5}:0.7) and (\B{7}:0.7) .. (\B{7}:1) ;
\draw[red, name=B11B13] (\B{11}:1) .. controls(\B{11}:0.7) and (\B{13}:0.7) .. (\B{13}:1) ;
\draw[red, name=B15B1] (\B{15}:1) .. controls(\B{15}:0.7) and (\B{1}:0.7) .. (\B{1}:1) ;
\end{tikzpicture}
\quad
 \begin{tikzpicture}[scale=1]
\newcommand{\degree}[0]{16}
\newcommand{\last}[0]{15}
\newcommand{\B}[1]{#1*360/\degree}
\draw (0,0) circle (1) ;
\draw[blue, name=B0B2] (\B{0}:1) .. controls(\B{0}:0.7) and (\B{2}:0.7) .. (\B{2}:1) ;
\draw[blue, name=B4B6] (\B{4}:1) .. controls(\B{4}:0.7) and (\B{6}:0.7) .. (\B{6}:1) ;
\draw[blue, name=B8B10] (\B{8}:1) .. controls(\B{8}:0.7) and (\B{10}:0.7) .. (\B{10}:1) ;
\draw[blue, name=B12B14] (\B{12}:1) .. controls(\B{12}:0.7) and (\B{14}:0.7) .. (\B{14}:1) ;
\draw[red, name=B1B3] (\B{1}:1) .. controls(\B{1}:0.7) and (\B{3}:0.7) .. (\B{3}:1) ;
\draw[red, name=B5B7](\B{5}:1) .. controls(\B{5}:0.7) and (\B{7}:0.7) .. (\B{7}:1) ;
\draw[red, name=B9B15] (\B{9}:1) .. controls(\B{9}:0.7) and (\B{15}:0.7) .. (\B{15}:1) ;
\draw[red, name=B11B13] (\B{11}:1) .. controls(\B{11}:0.7) and (\B{13}:0.7) .. (\B{13}:1) ;
\end{tikzpicture}
\quad
 \begin{tikzpicture}[scale=1]
\newcommand{\degree}[0]{16}
\newcommand{\last}[0]{15}
\newcommand{\B}[1]{#1*360/\degree}
\draw (0,0) circle (1) ;
\draw[blue, name=B0B2] (\B{0}:1) .. controls(\B{0}:0.7) and (\B{2}:0.7) .. (\B{2}:1) ;
\draw[blue, name=B4B6] (\B{4}:1) .. controls(\B{4}:0.7) and (\B{6}:0.7) .. (\B{6}:1) ;
\draw[blue, name=B8B14] (\B{8}:1) .. controls(\B{8}:0.7) and (\B{14}:0.7) .. (\B{14}:1) ;
\draw[blue, name=B10B12] (\B{10}:1) .. controls(\B{10}:0.7) and (\B{12}:0.7) .. (\B{12}:1) ;
\draw[red, name=B1B3] (\B{1}:1) .. controls(\B{1}:0.7) and (\B{3}:0.7) .. (\B{3}:1) ;
\draw[red, name=B5B15](\B{5}:1) .. controls(\B{5}:0.7) and (\B{15}:0.7) .. (\B{15}:1) ;
\draw[red, name=B7B9] (\B{7}:1) .. controls(\B{7}:0.7) and (\B{9}:0.7) .. (\B{9}:1) ;
\draw[red, name=B11B13] (\B{11}:1) .. controls(\B{11}:0.7) and (\B{13}:0.7) .. (\B{13}:1) ;
\end{tikzpicture}

\vspace{10pt}

\begin{tikzpicture}[scale=1]
\newcommand{\degree}[0]{16}
\newcommand{\last}[0]{15}
\newcommand{\B}[1]{#1*360/\degree}
\draw (0,0) circle (1) ;
\draw[blue, name=B0B2] (\B{0}:1) .. controls(\B{0}:0.7) and (\B{2}:0.7) .. (\B{2}:1) ;
\draw[blue, name=B4B6] (\B{4}:1) .. controls(\B{4}:0.7) and (\B{6}:0.7) .. (\B{6}:1) ;
\draw[blue, name=B8B14] (\B{8}:1) .. controls(\B{8}:0.7) and (\B{14}:0.7) .. (\B{14}:1) ;
\draw[blue, name=B10B12] (\B{10}:1) .. controls(\B{10}:0.7) and (\B{12}:0.7) .. (\B{12}:1) ;
\draw[red, name=B1B3] (\B{1}:1) .. controls(\B{1}:0.7) and (\B{3}:0.7) .. (\B{3}:1) ;
\draw[red, name=B5B7] (\B{5}:1) .. controls(\B{5}:0.7) and (\B{7}:0.7) .. (\B{7}:1) ;
\draw[red, name=B9B15] (\B{9}:1) .. controls(\B{9}:0.7) and (\B{15}:0.7) .. (\B{15}:1) ;
\draw[red, name=B11B13] (\B{11}:1) .. controls(\B{11}:0.7) and (\B{13}:0.7) .. (\B{13}:1) ;
\end{tikzpicture}
\quad
 \begin{tikzpicture}[scale=1]
\newcommand{\degree}[0]{16}
\newcommand{\last}[0]{15}
\newcommand{\B}[1]{#1*360/\degree}
\draw (0,0) circle (1) ;
\draw[blue, name=B0B2] (\B{0}:1) .. controls(\B{0}:0.7) and (\B{2}:0.7) .. (\B{2}:1) ;
\draw[blue, name=B4B6] (\B{4}:1) .. controls(\B{4}:0.7) and (\B{6}:0.7) .. (\B{6}:1) ;
\draw[blue, name=B8B10] (\B{8}:1) .. controls(\B{8}:0.7) and (\B{10}:0.7) .. (\B{10}:1) ;
\draw[blue, name=B12B14] (\B{12}:1) .. controls(\B{12}:0.7) and (\B{14}:0.7) .. (\B{14}:1) ;
\draw[red, name=B1B3] (\B{1}:1) .. controls(\B{1}:0.7) and (\B{3}:0.7) .. (\B{3}:1) ;
\draw[red, name=B5B15](\B{5}:1) .. controls(\B{5}:0.7) and (\B{15}:0.7) .. (\B{15}:1) ;
\draw[red, name=B7B9] (\B{7}:1) .. controls(\B{7}:0.7) and (\B{9}:0.7) .. (\B{9}:1) ;
\draw[red, name=B11B13] (\B{11}:1) .. controls(\B{11}:0.7) and (\B{13}:0.7) .. (\B{13}:1) ;
\end{tikzpicture}
\quad
 \begin{tikzpicture}[scale=1]
\newcommand{\degree}[0]{16}
\newcommand{\last}[0]{15}
\newcommand{\B}[1]{#1*360/\degree}
\draw (0,0) circle (1) ;
\draw[blue, name=B0B2] (\B{0}:1) .. controls(\B{0}:0.7) and (\B{2}:0.7) .. (\B{2}:1) ;
\draw[blue, name=B4B6] (\B{4}:1) .. controls(\B{4}:0.7) and (\B{6}:0.7) .. (\B{6}:1) ;
\draw[blue, name=B8B14] (\B{8}:1) .. controls(\B{8}:0.7) and (\B{14}:0.7) .. (\B{14}:1) ;
\draw[blue, name=B10B12] (\B{10}:1) .. controls(\B{10}:0.7) and (\B{12}:0.7) .. (\B{12}:1) ;
\draw[red, name=B1B7] (\B{1}:1) .. controls(\B{1}:0.7) and (\B{7}:0.7) .. (\B{7}:1) ;
\draw[red, name=B9B15](\B{9}:1) .. controls(\B{9}:0.7) and (\B{15}:0.7) .. (\B{15}:1) ;
\draw[red, name=B3B5] (\B{3}:1) .. controls(\B{3}:0.7) and (\B{5}:0.7) .. (\B{5}:1) ;
\draw[red, name=B11B13] (\B{11}:1) .. controls(\B{11}:0.7) and (\B{13}:0.7) .. (\B{13}:1) ;
\end{tikzpicture}
\end{itemize}

\subsection{Combinatorial closure of a signature}

\medskip

The existence of the notion of incidence relation between strata, which leads to the notion of poset and chain. 
A poset (partially ordered set) is a set $P$ together with a binary relation $\prec$. The elements $x$ and $y$ are comparable if $x \prec y$ and/or $y \prec x$ hold. A chain in a poset $P$ is a subset $C \subseteq P$ such that any two elements in $C$ are comparable.

\smallskip

We write $\sigma \prec \tau$, when $\tau$ can be obtained from $\sigma$ by a sequence of repeated contracting Whitehead moves. The signature $\tau$ is incident to $\sigma$. 
\begin{defn}
We call the union $\overline{\sigma}=\{\tau : \sigma \prec \tau\}$ the combinatorial closure of $\sigma$, and we define \[A_{\overline{\sigma}}=\cup_{\sigma \prec \tau} A_{\tau}.\] 
\end{defn}

\smallskip

\begin{defn}
Let $S$ be the set of signatures.
An upper bound in a subset of $(S,\prec)$ is a signature such that there exists no other signature $\tau$ in this subset verifying $\sigma\prec \tau$.  
 \end{defn}

\smallskip

\begin{lem}\label{L:123}
\vskip .1cm
 Let $\tau$ be a signature with a given  
 intersection of two red (or blue) diagonals. There are exactly two ways to smooth the intersection, 
which give two non-isotopic signatures $\sigma_1$ and $\sigma_2$ such that 
$\tau$ is incident to both $\sigma_1$ and $\sigma_2$. 
\vskip .1cm
\end{lem}

\begin{proof}
The intersection of two diagonals gives a double point and the argument to prove this follows from the desingularisation of double points.
Let $z_0$ be a double point in the drawing (or signature). The desingularisation consists of the topological model at $z_0$ consists in considering a suitable small open neighbourhood $U$ of $z_0$ and substituting $T\cap U$ for two disjoint open curves in such a way that one obtains a drawing with a double point less. This operation goes by the name of ``flip'' in and \cite{GCR,GCS}. There are exactly two ways of desingularising a double point in a real algebraic curve. The reason why the signatures are different is that the 4 asymptotic directions of the pair of curves are paired in the opposite way. 
\end{proof}

Similarly, we have the following for contracting Whitehead moves.  
\begin{lem}~\label{L:countWhitehead}
 Let $\tau_1$ and $\tau_2$ be obtained from
a signature $\sigma$ by two different complete contracting Whitehead moves.  Then,
$\tau_1$ and $\tau_2$ are different. 
\end{lem}
\begin{proof}
Let $\tau_1$ and $\tau_2$ be obtained from a signature $\sigma$ from two different contracting Whitehead moves. Two contracting Whitehead moves are different if they operate on different sets of edges. Let us consider $m_1$ (resp. $m_2$) edges of $\sigma$ lying in the boundary of cell of $\C$. Suppose that their set of leaves is $\{i_1,...i_{2m_{1}}\}$ (resp. $\{j_1,...j_{2m_{2}}\}$). A contracting Whitehead move glues those edges at one point. This gives a (star shaped) tree with leaves in the set $\{i_1,...i_{2m}\}$ (resp.$\{j_1,...j_{2m_{2}}\}$).  This gives the signature $\tau_1$ (resp.  $\tau_2$). Clearly $\tau_1$ can not be isotopic to $\tau_2$, with respect to the leaves, since  $\{i_1,...i_{2m}\}$ is different from $\{j_1,...j_{2m_{2}}\}$.
\end{proof}

\smallskip 

\begin{lem}\label{L:tree}
Let $\tau$ be a non generic signature. Consider in $\tau$ a vertex, incident to $m>2$ red 
(or blue)  curves. The signatures, obtained from $\tau$ in one single smoothing Whitehead move, are all distinct. 
\end{lem}
\begin{proof}
Locally, in a small neighbourhood around the intersection, the graph resembles a star shaped tree with one inner node and $2m$ leaves. Label those leaves $1,..,2m$. 
After one edge smoothing partial Whitehead move, the graph remains connected (one has simply detached a given branch and slided it along another one creating an extra nodal point. The valency of the first node being shrinked to $2n-2$). The topology of the graph is different. It is clear that  for different splittings the graphs are not isotopic. 
After a complete smoothing Whitehead move, the graph is disconnected: there exists a star-like tree  and a tree with at least one edge. There are $2m$ ways of inserting an edge if the two newly created vertices have different degrees and $m$ ways if they have same degree. These graphs are non-isotopic, with respect to the asymptotic directions of the leaves $2m$. 
\end{proof}

 \begin{lem}\label{L:123}
\vskip .1cm
Let $\tau$ be a signature of codimension $k$ with a given red 
(or blue) intersection of $m>2$ curves. Then there exists $Cat(m)$ distinct signatures which are smoothings of $\tau$ and of codimension $(k-(2m-3))$ obtained by complete smoothing Whitehead moves, where $Cat(m)$ is the $m^{th}$-Catalan number. 
\end{lem}
\begin{proof}
Draw in the neighborhood of the critical point $p$ a $2m$-gon. The vertices are the boundaries of the $m$ curves. Since we only consider the combinatorics, we may assume that the $2m$-gon is regular. Ungluing those $m$ curves gives $m$ disjoint curves in the regular $2m$-gon and the number of such possibilities is given by the $m^{th}$-Catalan number (see~\cite{Sta01}). 
\end{proof}

\subsection{Topological closure of a stratum}
\smallskip

\noindent \begin{lem}~\label{L:ball}
Let $\tau$ be a non-generic signature and let
$A_\tau$ be the corresponding strata of $\dpol_{d}$.  Let $P_0\in A_\tau$.
Then for every generic signature $\sigma$ such that $\tau$ is incident to 
$\sigma$, every $2d$-dimensional open ball $B_{\epsilon}$ containing $P_0$ intersects the
generic stratum $A_\sigma$.
\end{lem}

\begin{proof}
Suppose that there exists a generic signature $\sigma$ which does not verify that $B_{\epsilon}$ containing $P_0$ intersects the
generic stratum $A_\sigma$. 
Then, by Lemma\, \cite{L:6} there does not exist any path from a point $y \in A_{\tau}$ to $x\in A_{\sigma}$. Therefore $A_{\sigma}$ and $A_{\tau}$ are disjoint. In particular this implies that  $A_{\tau}$ is not in the closure of  $A_{\sigma}$. So, $\tau$ is not incident to $\sigma$.
\end{proof}

\begin{lem}
Let $\sigma$ be a generic signature and let
$\tau\ne\sigma$.  Then $\tau$ is incident to $\sigma$ if and only if the
following holds: for every pair of points $x,y\in \dpol_{d}$ with
$x\in A_\sigma$ and $y\in A_\tau$, there exists a continuous path 
$\gamma:[0,1]\rightarrow \dpol_n$ such that $\gamma(0)=x$, $\gamma(1)=y$ and 
$\gamma(t)\in A_\sigma$ for all $t\in [0,1)$. Any other such path $\rho$ 
from a point $x'\in A_\sigma$ to a point $y'\in A_\tau$ is 
homotopic to $\gamma$.
\end{lem}

\begin{proof}
 The result follows from lemma~\ref{L:ball}.  Indeed, 
if $\tau$ is incident to $\sigma$ then, a ball around $y\in
A_\tau$ necessarily intersects $A_\sigma$ and therefore there is a 
path from $y$ to a point $x'\in A_\sigma$. Composing this with a path
from $x'$ to $x$ in $A_\sigma$ we obtain $\gamma$.  Conversely,
if there is a path $\gamma$ from $x\in A_\sigma$ to $y\in A_\tau$,
then it is impossible to have an open ball containing
$y$ that does not intersect $A_\sigma$.
\end{proof}

We return to the original polynomials in order to define Whitehead moves in an analytic and then topological fashion.

\begin{prop}\label{D:sigmabar}
If $\sigma$ is a signature, 
the closure $\overline{A}_\sigma$ of
the stratum $A_\sigma$ in $\dpol_n$ is given by $$\overline{A}_\sigma=\cup_{\tau\in \overline{\sigma}} A_\tau,$$ where the {\it boundary}  of
$\sigma$ denoted $\overline{\sigma}$ consists of all incident signatures $\tau$. 
\end{prop}

\begin{proof}

One direction is easy. Indeed, if $x\in A_{\tau}$ where $\tau$ is incident to $\sigma$ then every $ 2d$-dimensional open set containing $x$ must intersect $A_{\sigma}$, so  $\cup_{\tau\in \overline{\sigma}} A_\tau\subset \overline{A}_\sigma$. 

For the other direction, let $x \in \overline{A}_\sigma\setminus A_{\sigma}$ and let $\tau$ be the signature of $x$. We first note that the dimension of $\tau$ can not be equal to the dimension of $\sigma$ because if they were equal, $A_{\tau}$ would be an open stratum disjoint from $A_{\sigma}$.

\smallskip

Therefore the dimension of $\tau$ is less than the dimension of $\sigma$. Let $U$ be any small open neighborhood of $x$. Let $y\in U\cap A_\sigma $ and let $\gamma$ be a path from $y$ to $x$ such that $\gamma \setminus x\subset A_{\sigma}$. Then every point $z \in \gamma\setminus x$ has the same signature $\sigma$. Using theorem~\ref{polyhedral_complex}, any path from the interior  to any point not in $A_{\sigma}$ must pass through the boundary of the polytope. 
Therefore any sequence of Whitehead moves and smoothings from $\sigma$ to $\tau$ must begin with Whitehead moves bringing $\sigma$ to a  signature which is incident to $\sigma$. But $x$ is the first point on  $\gamma $ where the signature changes and therefore $\tau$ must be incident to $\sigma$.  
\end{proof}

\smallskip

\smallskip

\begin{thm}~\label{T:cell}
The set of all $n-$signatures determine a stratification of the configuration space $Conf_n(\C)$ in which a signature of dimension $k$ corresponds to a stratum of dimension $k$ in $Conf_n(\C)$. 
\end{thm}
\begin{proof}
A signature corresponds uniquely to a class of drawings, a
class of drawings describes a region of polynomials (i.e. the set of polynomials having
the same signature) and the statement from theorem~\ref{T:Contract} shows that these regions are contractible. Due to proposition \ref{D:sigmabar}, lemma \ref{L:ball} and lemma \ref{L:ret}, we have a filtration verifying the definition of a stratified space. \
\end{proof}

\medskip

\begin{thm}
Let $\sigma$ be a generic signature.
Then, $\overline{A}_{\sigma}$ is a Goresky--MacPherson stratified space. 
\end{thm}
\begin{proof}
We have the following filtration: $\overline{A}_{\tau_0}\subset \overline{A}_{\tau_{1}}\subset \dots \subset \overline{A}_{\tau_{n-1}} \subset\overline{A}_{\tau_{n}}=\overline{A}_{\sigma}$. Take a point $p$ in  $\overline{A}_{\tau_{i}}\setminus \overline{A}_{\tau_{i-1}}$. Then, by Lem \ref{L:DeR}, $A_\sigma$ retracts onto $\overline{A}_{\tau_{i}}$ and we have a topological cone structure $cone(L)$, by the construction presented in the proof of Lem.~\ref{L:smo} where $L$ is a compact Hausdorff space. From the construction of our stratification, it follows that $L$ is endowed with an $n-i-1$ topological stratification: 
\[ L=L_{n-i-1}\supset \dots \supset L_{1}\supset L_{0} \supset L_{-1}=\emptyset.\] 

Now to show that locally a small neighborhood of a point $P$ in $\overline{A}_{\tau_{i}}$ behaves as an $i$-dimensional Euclidean space it is enough to consider the ramified cover $\pi_w:\dpol_n\to \tilde{\C}_{n}$. In particular, for the case where we are not on the ramification locus, we consider  an $(n-1)$-tuple of critical values $(w_1,\cdots w_{n-1})$ corresponding to the polynomial $P$ and take a small neighbourhood $U_w$ around it in the affine space $\C^{n-1}_{w}$. Then, by definition the pullback under $\pi_w$ over this neighbourhood is isomorphic to a product bundle $U_w\times \pi_w^{-1}(w)$. More generally, we have that points in the stratum cover a subset of an affine space and multiplicities correspond to proper subspaces of the affine space.

 So, there exists a distinguished neighborhood $N$ of $P$ in $\overline{A}_{\sigma}$ and a homeomorphism $\phi$ such that:
\[\phi:\Rl^i\times cone^{\circ}(L)\to N,\]  ($cone^{\circ}(L_{j})$ denotes the open cone) which takes each $\Rl^i\times cone^{\circ}(L_{j})$ homeomorphically to $N\cap \overline{A}_{\tau_{i+j+1}}$.

\end{proof}
\subsection{Multiple intersections of closures of strata}
This section investigates a combinatorial method to study multiple intersections between closures of strata indexed by generic signatures. This shows that a finite non-empty intersection of closures of generic strata has a lowest upper bound. This part is an independent investigation of the properties of the stratification, and is not directly used for the construction of the \v Cech cover. 

\smallskip

For simplicity, instead of signatures, we will use diagrams i.e.  embedded forests in the complex plane, such that their leaves lie on the boundary of a disc, on the $4n$ roots of the unity. We label the leaves in the trigonometric sense. We focus essentially on the combinatorics. So, from now on, we can assume that all the discs are of the same canonical radius and that their leaves coincide on the $4n$ roots of the unity.

\smallskip

The core idea of the construction is to superimpose generic diagrams $\sigma_0,...,\sigma_p$ such that leaves $\sqrt{1}^{4n}$ (and their labels) coincide. Then, we apply contracting Whitehead moves only to those diagonals, of the generic signatures, which in the superimposition are non identical.  
 
\smallskip

\subsubsection{Admissible superimposition of diagrams}

 \smallskip
 
Let $\sigma_0,\dots, \sigma_p$ be generic signatures and let $\Theta$ denote their superimposition. This superimposition is not well-defined as the diagonals of the different signatures can be positioned differently with respect to each other, since the diagonals of each signature are given only up to isotopy, but we will consider only those having the following properties:
 \begin{enumerate}
 \item all intersections are crossings (but not tangents) of at most two diagonals,
\item  the superimposition $\Theta$ cuts the disk into polygonal regions; we require that no region is a bigon,
\item the number of crossings of a given diagonal with other diagonals of  $\Theta$ must be minimal, with respect to a possible isotopy. We take representatives of isotopy classes of arcs. 
\end{enumerate}

\vskip.1cm 

We call such superimpositions {\it admissible}. 

\begin{lem}
Let $\sigma_0,\dots, \sigma_p$ be generic signatures. If there exists no superimposition $\Theta$ with the property that one diagonal has more than $p+1$ intersections with diagonals of the opposite color then there is no signature $\tau$ incident to all the $\sigma_{i}$.
 \end{lem}
 \begin{proof} The key point is the following. If $\sigma_{0},...,\sigma_{p}$ admit a signature $\tau$ incident to all of them, then $\tau$ has the following property: every segment of the tree $\tau$ (a segment is the part of an edge contained between two vertices, including leaves) belongs to at most one diagonal $(i,j)$ of each $\sigma_i$.  Thus, in particular, each segment can be considered as belonging to at most $p+1$ diagonals, one from each $\sigma_i$.  Thus, if
 a red diagonal of $\tau$ crosses $p+2$ or more blue diagonals in the superimposition, there is no one segment of $\tau$ which can belong to all of them, so
 the red diagonal will necessarily cross more than one blue segment of $\tau$, which is impossible.  \end{proof}

We will say that the set of signatures $\sigma_0,\dots, \sigma_p$ is {\it compatible} if it admits an admissible superimposition $\Theta$  with the property that no diagonal crosses more than $p+1$ diagonals of the opposite color. Note that a red diagonal  can never cross a blue diagonal in more than one point. Compatible sets of generic signatures may potentially have non-empty intersection.
We will now show how to give a condition on $\Theta$ to see whether or not this is the case.

\subsubsection{Graph associated to an intersection of generic signatures.} 
 
\medskip
 
 Let $\sigma_0,\dots, \sigma_p$ be a set of compatible generic signatures and let $\Theta$ be an  admissible superimposition.
Then $\Theta$ cuts the disk into polygonal regions.  Color a region red (resp. blue) if all its edges are red (resp. blue); the intersecting regions are purple.  Construct a graph from $\Theta$ as follows : 
place a vertex in each red or blue region (but not purple) with number of sides greater than three. 
 If two vertices lie in blue (resp. red) polygons that meet at a point, join them with a blue (resp. red) edge (even if this edge crosses purple regions). If two vertices lie in blue (resp. red) polygons that intersect along an edge of the opposite color, connect them with a blue (resp. red) edge. If two vertices lying in the same red (resp. blue) polygon can be connected by a segment inside the polygon which crosses only one purple region, add this segment. Connect each vertex to every terminal vertex lying in the same red (resp. blue) region, and also to any terminal vertex which can be reached by staying within the original red (resp. blue) polygon but crossing through a purple region formed by two blue (resp. red) diagonals emerging from that terminal vertex.  Finally, if any vertex of the graph has valency 2, we ignore this vertex and consider the two emerging edges as forming a single edge. We call this graph the graph associated to the superimposition.  
 
  \smallskip
 
 \begin{lem}
 The graph associated to $\Theta$ is independent of the actual choice of admissible superimposition $\Theta$.
 \end{lem}
 \begin{proof}
Let $\Theta$  and $\Theta'$ be  admissible superimpositions, and consider a given diagonal $D$. By the admissibility conditions the number of crossings of $D$ with diagonals of the other color is equal in $\Theta$  and $\Theta'$, and in fact the set of diagonals of the other color crossed by $D$ is identical in $\Theta$  and $\Theta'$. Therefore, the only possible modifications of the $\Theta$ is to move $D$ across an intersection of two diagonals of the other color. But this does not change the associated graph. 
 \end{proof}
 
\subsubsection{Compatible signatures}

\medskip

 \begin{defn} 
 The canonical graph associated to a set of compatible signatures $\sigma_0,\dots,\sigma_p$ is the graph associated to any  admissible superimposition $\Theta$. 
 \end{defn}

 \begin{thm}\label{Th:com}
Let $\sigma_0,\dots, \sigma_p$ be generic signatures. Then, there exists a signature incident to all the $\sigma_i$ if and only if the set $\sigma_0,\dots, \sigma_p$ is compatible and the associated canonical graph is a signature.
\end{thm}

\begin{proof} 
Replacing a blue (or red) polygon by a graph having the shape of a star as in the construction above involves diagonals of the different $\sigma_i$ which must be identified if we want to construct a common incident signature. The contracting moves may be stronger than strictly necessary (i.e. the signature $\tau$ may not be the signature of maximal dimension in the intersection), but any signature in the intersection must either have the same connected components as $\tau$, i.e. be obtained from $\tau$ by applying only smoothing Whitehead moves which do not increase the number of connected components of $\tau$ (these smoothings are the {\it partial smoothings}), or lie in the closures of these. Thus, up to such smoothings, the moves constructing $\tau$ are necessary in order to identify the diagonals of the $\sigma_i$.   

The contracting moves in the construction of the graph associated to $\Theta$, restricted to just one of the signatures $\sigma_i$, has the effect of making a contracting Whitehead move on the blue (resp. red) curves of this signature.  Thus, on each of the signatures, the graph construction reduces to a sequence of contracting Whitehead moves.  Thus the $\sigma_i$ possess a common incident signature if and only if $\tau$ is such a signature.
\end{proof}

\begin{rem}
In the appendix we show the classification of intersections between polygons which are allowed in order to have a non-empty intersection between closures of strata. At the end we show how to construct the canonical graph from the superimposition of a red and a blue polygon. 
\end{rem}

\medskip

\section{Good \v Cech cover}
This section shows how to construct the \v Cech cover, and includes a proof that multiple intersections are contractible.

\smallskip

We will introduce the notation $\underline{A_{\sigma}}$, for the set of all  elements $A_{\tau}$ in the set of strata verifying $\tau\prec \sigma$ i.e. 
$\underline{A_{\sigma}}=\cup_{\tau\prec\sigma} A_{\tau}$ and $codim(\tau)<codim(\sigma)$. The set of elements verifying $\tau\prec\sigma$, where  $codim(\tau)<codim(\sigma)$ is denoted by $\underline{\sigma}$. 

\medskip

\subsection{Contractibility theorem}
\begin{thm}~\label{T:clostra}
Let $\sigma$ be a non-generic signature . 
Then, $\underline{A_{\sigma}}$ is a contractible set.
\end{thm}
\begin{proof} 
We construct a diagram of spaces, where we use the simplicial realization of the poset structure $\underline{A_{\sigma}}$ and the contractibility arguments on the strata. 

For clarity, we discuss an inductive type of reasoning on a chain in the poset. 
An element $\tau$ lying in $\underline{A_{\sigma}}$ verifies  $\tau \prec \sigma$. In other words, there exists one unique  contracting Whitehead move from $\tau$ to $\sigma$ (uniqueness follows from lemma \ref{L:countWhitehead} and lemma \ref{L:tree}).
Suppose that $ codim(\sigma)=codim(\tau)+1$. Now by the deformation retract lemma \ref{L:DeR}, $A_{\sigma}$ is a deformation retract of $A_{\tau}$.

\smallskip 

Let $\tau'$ in  $\underline{A_{\sigma}}$ be of codimension $k$ and suppose that $A_{\sigma}$ is a deformation retract of $A_{\tau'}$. Let us show that for $ \tau''\prec \sigma$, where $codim(A_{\tau''})=codim(A_{\tau'})+1$, $A_{\sigma}$ is a deformation retract of $A_{\tau''}$. Now, by  the deformation retract lemma \ref{L:DeR}, $A_{\tau'}$ is a deformation retract of $A_{\tau''}$. So, by the induction hypothesis, we have that $A_{\sigma}$ is a deformation retract of $A_{\tau''}$. 

\smallskip 

The argument of contractibility concerning  $\underline{A_{\sigma}}$ comes from the simplicial realisation of this poset. In particular, we invoke the fact that there exists a unique maximal element in this poset (namely $\sigma$), giving it a cone structure. Now,  $A_{\tau}$ is contractible for each $\tau \preceq \sigma$ and whenever $\tau'\prec\tau$ we have that the induced map $A_{\tau}\to A_{\tau'}$ is a homotopy equivalence.
So, we can deduce that  $\underline{A_{\sigma}}$ is contractible. 
 \end{proof}

This proof easily adapts to give the following result. 

\begin{cor}\label{C:contra}
Let $\sigma$ be any signature and let $\tilde{\sigma}$ be a subset of $\underline{\sigma}$ with a maximal (or minimal) element.
Let $A_{\tilde{\sigma}}$ be  the union of $A_{\rho}$ with $\rho$ in $\tilde{\sigma}$. Then, $A_{\tilde{\sigma}}$ is contractible. 
\end{cor}

\medskip

\subsection{Thickening of strata}
Let $(\dpol_n, \mathcal{S})$ be the stratified smooth topological space. A system of  tubes is a family $T$ of triplets $(T_\sigma,\pi_\sigma,\rho_\sigma)_{\sigma\in S}$, where 
\begin{enumerate}
\item $T_\sigma$ is an open neighbourhood of the stratum $A_{\sigma}$, called its  tubular neighbourhood;

\item $\pi_\sigma$ is a strong deformation retraction $T_\sigma \to A_{\sigma}$ (i.e. a homotopy $h_t:T_{\sigma}\times [0,1]\to T_{\sigma}$ such that $h_0$ is the identity $h_1(T_{\sigma})=A_{\sigma}$ and  $h_t\arrowvert_{A_{\sigma}}=Id$ for all $t\in [0,1].$

\item$\rho_\sigma$ is a continuous function $T_\sigma \to \Rl_{\geq}$, called the distance function of the stratum, such that $\overline{A}_{\sigma} = \rho^{-1}(0)$;

\item for any pair of indices $\tau\prec \sigma$ the restriction of  $(\pi_\tau,\rho_\tau)$ to 
$T_\sigma \cap A_{\tau} \to A_{\sigma} \times \Rl_{\geq}$ is a smooth  submersion  

\item for all $x\in T_{\sigma} \cap T_\tau$ we have $\pi_{\sigma}(x)\in T_\tau$ ,$\pi_\tau\pi_\sigma(x)=\pi_\tau(x)$ and 
 $\rho_\tau(\pi_\sigma(x)) = \rho_\tau(x)$.
 \end{enumerate}

\medskip

\subsection{The \v Cech cover}
We consider only the red (resp. blue) curves of a signature $\sigma$ (i.e. the topological model). This is denoted by $\sigma^{R}$ (resp. $\sigma^{B}$).
\begin{defn}
A connected component in $\sigma^{R}$ (resp. $\sigma^B$) is called a tree.
\end{defn}

Consider the set of leaves of a signature and label them $\{1,...,4n\}$. Then the labels of a connected component (i.e. a tree) of a topological model form a subset (a block) of $\{1,...,4n\}$ and it is a set of even cardinality. A complete Whitehead move partitions this subset into two disjoint blocks of even cardinality each. These partitions are non crossing partitions. 

Two blocks of a partition are said to cross if there are integers $i < j < k < l$ such that $i,k$ belong to one block and $j, l$ belong to the other block. If no two blocks cross, then the partition is said to be noncrossing. A noncrossing matching (of order $n$) is a noncrossing partition (of a set of size $2n$) into blocks of size 2. This occurs in the generic case i.e. when there exist no real (or imaginary) critical values.   

\smallskip

\begin{lem}\label{L:con}

Consider two non-generic signatures $\sigma$ and $\tau$. Let $\mathcal{P}_1$ be the non crossing partition of blocks of leaves which form the connected components in $\sigma^{R}$ (resp. $\sigma^B$). Let $\mathcal{P}_2$ be the non crossing partition of blocks of leaves which form the connected components in $\tau^{R}$ (resp. $\tau^B$).
Then, $\underline{\sigma}\cap \underline{\tau}\ne \emptyset$ if the intersections of  $\mathcal{P}_1$ and $\mathcal{P}_2$ form a non crossing partition of the set of leaves and that each block has even cardinality.

\end{lem}
\begin{proof}
Suppose that the intersection of $\mathcal{P}_1$ and $\mathcal{P}_2$ does not form a crossing partition. Then there exists at least one block of  $\mathcal{P}_1$ and a block of $\mathcal{P}_2$ forming a crossing. But by construction two trees with a different set of leaves can never be modified, using smoothing Whitehead moves, into becoming identical. So, the non crossing condition is necessary. 
However this is not a sufficient condition, since the intersection of  $\mathcal{P}_1$ and $\mathcal{P}_2$ can give singletons and blocks of odd cardianlity. Now, this situation cannot make sense since we are considering connected components which by definition have a set of leaves of even cardinality.

\end{proof}

We proceed by the following algorithm.

\medskip

\begin{lem}\label{L:un}
Let $\sigma$ and $\tau$ be two distinct non-generic signatures. Consider $\underline{\sigma}:=\{\mu\, :\, \mu\prec \sigma\}$.
Suppose that $\underline{\sigma}\cap \underline{\tau}\ne \emptyset$. If the lowest upper bound exists in $\underline{\sigma}\cap \underline{\tau}$ then it is unique.  
\end{lem}
\begin{proof}
Consider the topological model $\sigma^{R}$. Label the leaves of this red graph by $\{1,\cdots 2n\}$.
To each connected component of $\sigma^{R}$ associate a block of leaves i.e. a subset of $\{1,\cdots, 2n\}$ defining the set of leaves of that tree. 
Then, this set of blocks forms a partition of $\{1,\cdots 2n\}$. If we denote each block by the set $A_i$ then we get the following non-crossing partition 
$\{1,\cdots, 2n\}=\sqcup_{i=1}^p A_i$ where $p$ is the number of blocks (that means of connected components in $\sigma^{R}$). Call this non-crossing partition $\mathcal{P}_1$.
Proceed the same way for $\tau^{R}$ and suppose that the set of blocks of non-crossing partitions is given by $\sqcup_{j=1}^r B_j$. Then there exist $r$ connected components in $\tau^{R}$ and the subsets of leaves of the $i$-th connected component lie in the given $i$-th block $B_i.$ Call this non-crossing partition $\mathcal{P}_2$.

Now, by hypothesis, we have that $\underline{\sigma}\cap \underline{\tau}\ne \emptyset$. So, regarding the relation between blocks in $\mathcal{P}_1$ and $\mathcal{P}_2$ they should have a non-crossing relation. For instance one block in $\mathcal{P}_1$ can be included in a bigger block of $\mathcal{P}_2$ (or vice versa) since this does not violate the non-crossing rule. We shall focus on the intersection of the blocks in $\mathcal{P}_1$ with those in $\mathcal{P}_2$.  This forms a new non-crossing partition $\mathcal{P}_3$ of the set of leaves.

Note that even though a pair of connected components in $\sigma^{R}$ and $\tau^{R}$ might have a coinciding set of leaves this does not provide a unique signature: to one block of leaves say $\{i_1,...,i_{2k}\}$ there corresponds a set of non isomorphic trees, if $k>1$. So, we remedy to this using the following procedure.

Let us consider a block in $\mathcal{P}_3$. This corresponds to a connected component with leaves in the set $\{i_1,...,i_{2k}\}$. Then, 
take the connected component $\sigma^R$ or $\tau^R$ having the minimal number of inner nodes (vertices which are not leaves!) and modify it, by applying a partial smoothing Whitehead move. Repeat until the number of vertices in both connected components coincide and in such a way that edge relations between leaves and inner nodes are identical in both connected components. Those trees are thus isomorphic and there exists one unique such graph. Apply exactly the same procedure for $\sigma^B$ and $\tau^B$ as previously. Therefore, by construction greatest lower bound is unique.
\end{proof}

\smallskip

We will use this construction to thicken all our strata. Concerning notation, a thickened stratum $A$  is denoted by $A^+$.
\begin{thm}
Let $\sigma_0,...,\sigma_p$ be a set of non-generic signature, being upper bounds.
Then, the open sets of the \v Cech cover are formed by the thickened sets  $\underline{A_{\sigma_i}}^{+}$.
\end{thm}

\begin{proof}
Indeed, to have a \v Cech cover, it is necessary to have open sets and that their multiple intersections are contractible. We have shown that $\underline{A_{\sigma}}$ is contractible. 
Applying the construction above, we thicken every stratum in $\underline{A_{\sigma}}$. Condition 4 of this construction implies that to cover our space we need the union of all thickened strata lying in  $\underline{A_{\sigma_i}}$, which we denote by $\underline{A_{\sigma}}^{+}$.

\smallskip 

Multiple intersections are contractible. 
For simplicity we consider the case of two intersecting sets, where $A=\underline{A_{\sigma}}^{+}$ and $B= \underline{A_{\sigma'}}^{+}$.
We have shown in the lemma \ref{L:un} that if the intersection is non-empty then there exists a unique greatest lower bound. 
So, we can identify $A\cap B$  to the tubular neighbourhood of $\underline{A_{\tau}}$, where $\tau$ is the meet of $\sigma$ and $\sigma'$ given by lemma~\ref{L:un}. Theorem \ref{T:clostra} states that $\underline{A_{\tau}}$ is contractible. The generalisation can be easily obtained by induction. Therefore, multiple intersections are contractible. 
\end{proof}

\vfill\eject

\appendix

\section{Superimposition of signatures}
Let $\sigma_0 \cup \dots \cup \sigma_p$ be compatible generic signatures and $\Theta$ denote an admissible superimposition.  In this subsection, we  digress briefly in order to give a visual description of the conditions on the superimposition $\Theta$ for the associated graph to be a signature. In fact, it is quite rare for signatures to intersect. Almost always the canonical graph will not be a signature. 
 Given a red polygon and a blue polygon of $\Theta$, they must either be disjoint or intersect in one of exactly four possible ways: 
\begin{itemize}
\item the intersection is a three sided polygon with two red (resp. blue) edges and one blue edge (resp. red); 
\begin{center}
\begin{tikzpicture}[scale=1.2]
\newcommand{\degree}[0]{6}
\newcommand{\last}[0]{11}
\newcommand{\B}[1]{#1*180/\degree}
\newcommand{\R}[1]{#1*180/\degree + 90/\degree}

\draw[blue, name=B4B9] (\B{4}:1) .. controls(\B{4}:0.1) and (\B{9}:0.7) .. (\B{9}:1) ;
\draw[red, name=R2R6] (\R{2}:1) .. controls(\R{2}:0.7) and (\R{6}:0.7) .. (\R{6}:1) ;
\draw[red, name=R6R11] (\R{6}:1) .. controls(\R{6}:0.7) and (\R{11}:0.7) .. (\R{11}:1) ;
\end{tikzpicture}
\end{center}

\item The intersection is a four sided polygon with two red and two blue edges; 
\begin{center}
\begin{tikzpicture}[scale=1.2]
\newcommand{\degree}[0]{6}
\newcommand{\last}[0]{11}
\newcommand{\XDOT}[0]{[black](i-1) circle (0.045)}
\newcommand{\B}[1]{#1*180/\degree}
\newcommand{\BDOT}[0]{[blue,fill=white](i-1) circle (0.05)}
\newcommand{\GDOT}[0]{[blue,fill=white] circle (0.15)}
\newcommand{\R}[1]{#1*180/\degree + 90/\degree}
\newcommand{\RDOT}[0]{[red,fill=white](i-1) circle (0.05)}
\newcommand{\HDOT}[0]{[red,fill=white] circle (0.15)}

\draw[blue, name=B0B1] (\B{0}:1) .. controls(\B{0}:0.7) and (\B{1}:0.7) .. (\B{1}:1) ;
\draw[blue, name=B0B11] (\B{0}:1) .. controls(\B{0}:0.7) and (\B{11}:0.7) .. (\B{11}:1) ;
\draw[blue, name=B1B5] (\B{1}:1) .. controls(\B{1}:0.7) and (\B{5}:0.7) .. (\B{5}:1) ;
\draw[blue, name=B7B11] (\B{7}:1) .. controls(\B{7}:0.7) and (\B{11}:0.7) .. (\B{11}:1) ;
\draw[blue, name=B5B6] (\B{5}:1) .. controls(\B{5}:0.7) and (\B{6}:0.7) .. (\B{6}:1) ;
\draw[blue, name=B6B7] (\B{6}:1) .. controls(\B{6}:0.7) and (\B{7}:0.7) .. (\B{7}:1) ;
\draw[red, name=R3R4] (\R{3}:1) .. controls(\R{3}:0.7) and (\R{4}:0.7) .. (\R{4}:1) ;
\draw[red, name=R2R3] (\R{2}:1) .. controls(\R{2}:0.7) and (\R{3}:0.7) .. (\R{3}:1) ;
\draw[red, name=R4R8] (\R{4}:1) .. controls(\R{4}:0.7) and (\R{8}:0.7) .. (\R{8}:1) ;
\draw[red, name=R2R10] (\R{2}:1) .. controls(\R{2}:0.7) and (\R{10}:0.7) .. (\R{10}:1) ;
\draw[red, name=R9R8] (\R{9}:1) .. controls(\R{9}:0.7) and (\R{8}:0.7) .. (\R{8}:1) ;
\draw[red, name=R9R10] (\R{9}:1) .. controls(\R{9}:0.7) and (\R{10}:0.7) .. (\R{10}:1) ;
\end{tikzpicture}
\end{center}
\item The intersection is two triangles joined at a point, formed by two crossing blue (resp. red) diagonals, cut transversally on either side of the intersection by two red (resp. blue) diagonals; note that this means that two polygons of the same color meet at a point and both intersect the polygon of the other color; \begin{center}
\begin{tikzpicture}[scale=1.2]
\newcommand{\degree}[0]{6}
\newcommand{\last}[0]{11}
\newcommand{\B}[1]{#1*180/\degree}
\newcommand{\R}[1]{#1*180/\degree + 90/\degree}
\draw[blue, name=B0B1] (\B{0}:1) .. controls(\B{0}:0.7) and (\B{1}:0.7) .. (\B{1}:1) ;
\draw[blue, name=B0B11] (\B{0}:1) .. controls(\B{0}:0.7) and (\B{11}:0.7) .. (\B{11}:1) ;
\draw[blue, name=B1B7] (\B{1}:1) .. controls(\B{1}:0.7) and (\B{7}:0.7) .. (\B{7}:1) ;
\draw[blue, name=B5B11] (\B{5}:1) .. controls(\B{5}:0.7) and (\B{11}:0.7) .. (\B{11}:1) ;
\draw[blue, name=B5B6] (\B{5}:1) .. controls(\B{5}:0.7) and (\B{6}:0.7) .. (\B{6}:1) ;
\draw[blue, name=B6B7] (\B{6}:1) .. controls(\B{6}:0.7) and (\B{7}:0.7) .. (\B{7}:1) ;
\draw[red, name=R3R4] (\R{3}:1) .. controls(\R{3}:0.7) and (\R{4}:0.7) .. (\R{4}:1) ;
\draw[red, name=R2R3] (\R{2}:1) .. controls(\R{2}:0.7) and (\R{3}:0.7) .. (\R{3}:1) ;
\draw[red, name=R4R8] (\R{4}:1) .. controls(\R{4}:0.7) and (\R{8}:0.7) .. (\R{8}:1) ;
\draw[red, name=R2R10] (\R{2}:1) .. controls(\R{2}:0.7) and (\R{10}:0.7) .. (\R{10}:1) ;
\draw[red, name=R9R8] (\R{9}:1) .. controls(\R{9}:0.7) and (\R{8}:0.7) .. (\R{8}:1) ;
\draw[red, name=R9R10] (\R{9}:1) .. controls(\R{9}:0.7) and (\R{10}:0.7) .. (\R{10}:1) ;
\end{tikzpicture}
\end{center}
\item The intersection is a point at which two blue diagonals and two red diagonals all cross in the cyclic order red, red, blue, blue, red, red, blue, blue; note that this means that in fact four polygons meet at a point, each being the same color as the opposing one.
\begin{center}
\begin{tikzpicture}[scale=1.2]
\newcommand{\degree}[0]{6}
\newcommand{\last}[0]{11}
\newcommand{\B}[1]{#1*180/\degree}
\newcommand{\R}[1]{#1*180/\degree + 90/\degree}
\draw[blue, name=B0B1] (\B{0}:1) .. controls(\B{0}:0.7) and (\B{1}:0.7) .. (\B{1}:1) ;
\draw[blue, name=B0B11] (\B{0}:1) .. controls(\B{0}:0.7) and (\B{11}:0.7) .. (\B{11}:1) ;
\draw[blue, name=B1B7] (\B{1}:1) .. controls(\B{1}:0.7) and (\B{7}:0.7) .. (\B{7}:1) ;
\draw[blue, name=B5B11] (\B{5}:1) .. controls(\B{5}:0.7) and (\B{11}:0.7) .. (\B{11}:1) ;
\draw[blue, name=B5B6] (\B{5}:1) .. controls(\B{5}:0.7) and (\B{6}:0.7) .. (\B{6}:1) ;
\draw[blue, name=B6B7] (\B{6}:1) .. controls(\B{6}:0.7) and (\B{7}:0.7) .. (\B{7}:1) ;
\draw[red, name=R3R4] (\R{3}:1) .. controls(\R{3}:0.7) and (\R{4}:0.7) .. (\R{4}:1) ;
\draw[red, name=R2R3] (\R{2}:1) .. controls(\R{2}:0.7) and (\R{3}:0.7) .. (\R{3}:1) ;
\draw[red, name=R2R8] (\R{2}:1) .. controls(\R{2}:0.7) and (\R{8}:0.7) .. (\R{8}:1) ;
\draw[red, name=R4R10] (\R{4}:1) .. controls(\R{4}:0.7) and (\R{10}:0.7) .. (\R{10}:1) ;
\draw[red, name=R9R8] (\R{9}:1) .. controls(\R{9}:0.7) and (\R{8}:0.7) .. (\R{8}:1) ;
\draw[red, name=R9R10] (\R{9}:1) .. controls(\R{9}:0.7) and (\R{10}:0.7) .. (\R{10}:1) ;
\end{tikzpicture}
\end{center}
\end{itemize}
From such a disposition of diagonals, we obtain the graph described above.
 The graph must be a forest with even valency at every non-terminal vertex.

\vspace{5pt}
In the following we show how to construct a canonical graph from the superimposition of a red and a blue polygon, using the second superimposition above.
\begin{figure}[ht]~\label{F:1}
\begin{center}
\begin{tikzpicture}[scale=0.9]
\node (a) at (-6,0) {
\begin{tikzpicture}[scale=1.2]
\newcommand{\degree}[0]{6}
\newcommand{\last}[0]{11}
\newcommand{\B}[1]{#1*180/\degree}
\newcommand{\R}[1]{#1*180/\degree + 90/\degree}
\draw[blue, name=B0B1] (\B{0}:1) .. controls(\B{0}:0.7) and (\B{1}:0.7) .. (\B{1}:1) ;
\draw[blue, name=B0B11] (\B{0}:1) .. controls(\B{0}:0.7) and (\B{11}:0.7) .. (\B{11}:1) ;
\draw[blue, name=B1B5] (\B{1}:1) .. controls(\B{1}:0.7) and (\B{5}:0.7) .. (\B{5}:1) ;
\draw[blue, name=B7B11] (\B{7}:1) .. controls(\B{7}:0.7) and (\B{11}:0.7) .. (\B{11}:1) ;
\draw[blue, name=B5B6] (\B{5}:1) .. controls(\B{5}:0.7) and (\B{6}:0.7) .. (\B{6}:1) ;
\draw[blue, name=B6B7] (\B{6}:1) .. controls(\B{6}:0.7) and (\B{7}:0.7) .. (\B{7}:1) ;
\draw[red, name=R3R4] (\R{3}:1) .. controls(\R{3}:0.7) and (\R{4}:0.7) .. (\R{4}:1) ;
\draw[red, name=R2R3] (\R{2}:1) .. controls(\R{2}:0.7) and (\R{3}:0.7) .. (\R{3}:1) ;
\draw[red, name=R4R8] (\R{4}:1) .. controls(\R{4}:0.7) and (\R{8}:0.7) .. (\R{8}:1) ;
\draw[red, name=R2R10] (\R{2}:1) .. controls(\R{2}:0.7) and (\R{10}:0.7) .. (\R{10}:1) ;
\draw[red, name=R9R8] (\R{9}:1) .. controls(\R{9}:0.7) and (\R{8}:0.7) .. (\R{8}:1) ;
\draw[red, name=R9R10] (\R{9}:1) .. controls(\R{9}:0.7) and (\R{10}:0.7) .. (\R{10}:1) ;
\end{tikzpicture}
};
\node (b) at (-2,0) {
\begin{tikzpicture}[scale=1.2]
\newcommand{\degree}[0]{6}
\newcommand{\last}[0]{11}
\newcommand{\B}[1]{#1*180/\degree}
\newcommand{\R}[1]{#1*180/\degree + 90/\degree}
\draw[blue, name=B0B1] (\B{0}:1) .. controls(\B{0}:0.7) and (\B{1}:0.7) .. (\B{1}:1) ;
\draw[blue, name=B0B11] (\B{0}:1) .. controls(\B{0}:0.7) and (\B{11}:0.7) .. (\B{11}:1) ;
\draw[blue, name=B1B5] (\B{1}:1) .. controls(\B{1}:0.7) and (\B{5}:0.7) .. (\B{5}:1) ;
\draw[blue, name=B7B11] (\B{7}:1) .. controls(\B{7}:0.7) and (\B{11}:0.7) .. (\B{11}:1) ;
\draw[blue, name=B5B6] (\B{5}:1) .. controls(\B{5}:0.7) and (\B{6}:0.7) .. (\B{6}:1) ;
\draw[blue, name=B6B7] (\B{6}:1) .. controls(\B{6}:0.7) and (\B{7}:0.7) .. (\B{7}:1) ;
\draw[red, name=R3R4] (\R{3}:1) .. controls(\R{3}:0.7) and (\R{4}:0.7) .. (\R{4}:1) ;
\draw[red, name=R2R3] (\R{2}:1) .. controls(\R{2}:0.7) and (\R{3}:0.7) .. (\R{3}:1) ;
\draw[red, name=R4R8] (\R{4}:1) .. controls(\R{4}:0.7) and (\R{8}:0.7) .. (\R{8}:1) ;
\draw[red, name=R2R10] (\R{2}:1) .. controls(\R{2}:0.7) and (\R{10}:0.7) .. (\R{10}:1) ;
\draw[red, name=R9R8] (\R{9}:1) .. controls(\R{9}:0.7) and (\R{8}:0.7) .. (\R{8}:1) ;
\draw[red, name=R9R10] (\R{9}:1) .. controls(\R{9}:0.7) and (\R{10}:0.7) .. (\R{10}:1) ;

\draw[blue, fill= blue] (-0.6,0) circle (0.05) ;
\draw[blue, fill= blue] (0.6,0) circle (0.05) ;
\draw[red, fill= red] (-0.15,0.58) circle (0.05) ;
\draw[red, fill= red] (0.1,-0.55) circle (0.05) ;

\end{tikzpicture}
};
\node (c) at (2,0) {
\begin{tikzpicture}[scale=1.2]
\newcommand{\degree}[0]{6}
\newcommand{\last}[0]{11}
\newcommand{\B}[1]{#1*180/\degree}
\newcommand{\R}[1]{#1*180/\degree + 90/\degree}
\draw[blue, name=B0B1] (\B{0}:1) .. controls(\B{0}:0.7) and (\B{1}:0.7) .. (\B{1}:1) ;
\draw[blue, name=B0B11] (\B{0}:1) .. controls(\B{0}:0.7) and (\B{11}:0.7) .. (\B{11}:1) ;
\draw[blue, name=B1B5] (\B{1}:1) .. controls(\B{1}:0.7) and (\B{5}:0.7) .. (\B{5}:1) ;
\draw[blue, name=B7B11] (\B{7}:1) .. controls(\B{7}:0.7) and (\B{11}:0.7) .. (\B{11}:1) ;
\draw[blue, name=B5B6] (\B{5}:1) .. controls(\B{5}:0.7) and (\B{6}:0.7) .. (\B{6}:1) ;
\draw[blue, name=B6B7] (\B{6}:1) .. controls(\B{6}:0.7) and (\B{7}:0.7) .. (\B{7}:1) ;
\draw[red, name=R3R4] (\R{3}:1) .. controls(\R{3}:0.7) and (\R{4}:0.7) .. (\R{4}:1) ;
\draw[red, name=R2R3] (\R{2}:1) .. controls(\R{2}:0.7) and (\R{3}:0.7) .. (\R{3}:1) ;
\draw[red, name=R4R8] (\R{4}:1) .. controls(\R{4}:0.7) and (\R{8}:0.7) .. (\R{8}:1) ;
\draw[red, name=R2R10] (\R{2}:1) .. controls(\R{2}:0.7) and (\R{10}:0.7) .. (\R{10}:1) ;
\draw[red, name=R9R8] (\R{9}:1) .. controls(\R{9}:0.7) and (\R{8}:0.7) .. (\R{8}:1) ;
\draw[red, name=R9R10] (\R{9}:1) .. controls(\R{9}:0.7) and (\R{10}:0.7) .. (\R{10}:1) ;

\draw[blue, name=B0B6] (\B{0}:1) .. controls(\B{0}:0.7) and (\B{6}:0.7) .. (\B{6}:1) ;
\draw[blue, name=B1B11] (\B{1}:1) .. controls(\B{1}:0.6) and (\B{11}:0.6) .. (\B{11}:1) ;
\draw[blue, name=B5B7] (\B{5}:1) .. controls(\B{5}:0.6) and (\B{7}:0.6) .. (\B{7}:1) ;

\draw[blue, fill= blue] (-0.6,0) circle (0.05) ;
\draw[blue, fill= blue] (0.6,0) circle (0.05) ;
\draw[red, fill= red] (-0.15,0.58) circle (0.05) ;
\draw[red, fill= red] (0.1,-0.55) circle (0.05) ;
\end{tikzpicture}
};
\node (d) at (6,0) {
\begin{tikzpicture}[scale=1.2]
\newcommand{\degree}[0]{6}
\newcommand{\last}[0]{11}
\newcommand{\B}[1]{#1*180/\degree}
\newcommand{\R}[1]{#1*180/\degree + 90/\degree}

\draw[blue, name=B0B1] (\B{0}:1) .. controls(\B{0}:0.7) and (\B{1}:0.7) .. (\B{1}:1) ;
\draw[blue, name=B0B11] (\B{0}:1) .. controls(\B{0}:0.7) and (\B{11}:0.7) .. (\B{11}:1) ;
\draw[blue, name=B1B5] (\B{1}:1) .. controls(\B{1}:0.7) and (\B{5}:0.7) .. (\B{5}:1) ;
\draw[blue, name=B7B11] (\B{7}:1) .. controls(\B{7}:0.7) and (\B{11}:0.7) .. (\B{11}:1) ;
\draw[blue, name=B5B6] (\B{5}:1) .. controls(\B{5}:0.7) and (\B{6}:0.7) .. (\B{6}:1) ;
\draw[blue, name=B6B7] (\B{6}:1) .. controls(\B{6}:0.7) and (\B{7}:0.7) .. (\B{7}:1) ;
\draw[red, name=R3R4] (\R{3}:1) .. controls(\R{3}:0.7) and (\R{4}:0.7) .. (\R{4}:1) ;
\draw[red, name=R2R3] (\R{2}:1) .. controls(\R{2}:0.7) and (\R{3}:0.7) .. (\R{3}:1) ;
\draw[red, name=R4R8] (\R{4}:1) .. controls(\R{4}:0.7) and (\R{8}:0.7) .. (\R{8}:1) ;
\draw[red, name=R2R10] (\R{2}:1) .. controls(\R{2}:0.7) and (\R{10}:0.7) .. (\R{10}:1) ;
\draw[red, name=R9R8] (\R{9}:1) .. controls(\R{9}:0.7) and (\R{8}:0.7) .. (\R{8}:1) ;
\draw[red, name=R9R10] (\R{9}:1) .. controls(\R{9}:0.7) and (\R{10}:0.7) .. (\R{10}:1) ;

\draw[blue, name=B0B6] (\B{0}:1) .. controls(\B{0}:0.7) and (\B{6}:0.7) .. (\B{6}:1) ;
\draw[blue, name=B1B11] (\B{1}:1) .. controls(\B{1}:0.6) and (\B{11}:0.6) .. (\B{11}:1) ;
\draw[blue, name=B5B7] (\B{5}:1) .. controls(\B{5}:0.6) and (\B{7}:0.6) .. (\B{7}:1) ;

\draw[red, name=R3R9] (\R{3}:1) .. controls(\R{3}:0.7) and (\R{9}:0.7) .. (\R{9}:1) ;
\draw[red, name=R2R4] (\R{2}:1) .. controls(\R{2}:0.6) and (\R{4}:0.6) .. (\R{4}:1) ;
\draw[red, name=R8R10] (\R{8}:1) .. controls(\R{8}:0.6) and (\R{10}:0.6) .. (\R{10}:1) ;
\draw[blue, fill= blue] (-0.6,0) circle (0.05) ;
\draw[blue, fill= blue] (0.6,0) circle (0.05) ;
\draw[red, fill= red] (-0.15,0.58) circle (0.05) ;
\draw[red, fill= red] (0.15,-0.58) circle (0.05) ;

\end{tikzpicture}
};
\draw[thick, ->] (a)--(b);
\draw[thick, ->] (b)--(c);
\draw[thick, ->] (c)--(d);
\end{tikzpicture}
\end{center}
\caption{Admissible superimposition}
\end{figure}
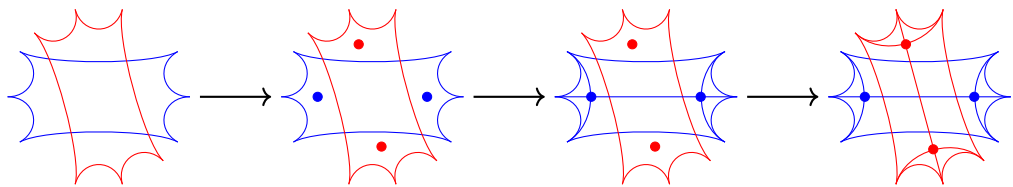
\begin{figure}[h]
\begin{center}
\begin{tikzpicture}[scale=1.2]
\newcommand{\degree}[0]{8}
\newcommand{\B}[1]{#1*180/\degree}
\newcommand{\R}[1]{#1*180/\degree }
 
\draw[blue, name=B0B8](\B{0}:1) .. controls(\B{0}:0.7)and (\B{8}:0.7) .. (\B{8}:1) ;
\draw[blue, name=R7R9](\R{7}:1) .. controls(\R{7}:0.3)and (\R{9}:0.3) .. (\R{9}:1) ;
\draw[blue, name=R15R1](\B{15}:1) .. controls(\B{15}:0.3)and (\B{1}:0.3) .. (\B{1}:1) ;

\draw[red, name=B12B4](\B{12}:1) .. controls(\B{12}:0.7)and (\B{4}:0.7) .. (\B{4}:1) ;
\draw[red, name=R3R5](\R{3}:1) .. controls(\R{3}:0.3)and (\R{5}:0.3) .. (\R{5}:1) ;
\draw[red, name=R11R13](\R{11}:1) .. controls(\R{11}:0.3)and (\R{13}:0.3) .. (\R{13}:1) ;

\end{tikzpicture}
\end{center}
\caption{Canonical graph}
\end{figure}
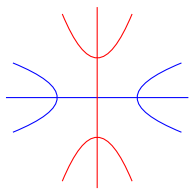

\vfill\eject

\end{document}